\documentclass{amsart}
\usepackage{amsmath}
\usepackage{graphicx}
\usepackage{amsfonts}
\usepackage{amssymb}

\newtheorem{theorem}{Theorem}
\theoremstyle{plain}

\newtheorem{definition}{Definition}
\newtheorem{example}{Example}

\newtheorem{proposition}{Proposition}
\newtheorem{remark}{Remark}

\numberwithin{equation}{section}
\input{tcilatex}

\begin{document}
\title{The Biderivative and $A_{\infty}$-Bialgebras}
\author{Samson Saneblidze$^{1}$}
\address{A. Razmadze Mathematical Institute\\
Georgian Academy of Sciences\\
M. Aleksidze st., 1\\
0193 Tbilisi, Georgia}
\email{sane@rmi.acnet.ge}
\author{Ronald Umble$^{2}$}
\address{Department of Mathematics\\
Millersville University of Pennsylvania\\
Millersville, PA. 17551}
\email{ron.umble@millersville.edu}
\thanks{$^{1} $ This research described in this publication was made
possible in part by Award No. GM1-2083 of the U.S. Civilian Research and
Development Foundation for the Independent States of the Former Soviet Union
(CRDF) and by Award No. 99-00817 of INTAS}
\thanks{$^{2}$ This research funded in part by a Millersville University
faculty research grant.}

\begin{abstract}
An $A_{\infty }$-bialgebra is a DGM $H$ equipped with structurally
compatible operations $\left\{ \omega ^{j,i}:H^{\otimes i}\rightarrow
H^{\otimes j}\right\} $ such that $\left( H,\omega ^{1,i}\right) $ is an $%
A_{\infty }$-algebra and $\left( H,\omega ^{j,1}\right) $ is an $A_{\infty }$%
-coalgebra. Structural compatibility is controlled by the biderivative
operator $Bd$, defined in terms of two kinds of cup products on certain
cochain algebras of pemutahedra over the universal PROP $U=End\left(
TH\right) $.
\end{abstract}

\date{March 19, 2005.}
\maketitle

\begin{center}
\textit{To Jim Stasheff on the occasion of his 68th birthday.}\vspace{0.2in}
\end{center}

\section{Introduction}

In his seminal papers of 1963, J. Stasheff \cite{Stasheff} introduced the
notion of an $A_{\infty}$-algebra, which is (roughly speaking) a DGA in
which the associative law holds up to homotopy. Since then, $A_{\infty}$%
-algebras have assumed their rightful place as fundamental structures in
algebra \cite{Markl}, \cite{Shnider}, topology \cite{Kadeishvili1}, \cite%
{Lin}, \cite{Stasheff2}, and mathematical physics \cite{Kimura1}, \cite%
{Kimura2}, \cite{Penkava1}, \cite{Penkava2}, \cite{Wiesbrock}, \cite%
{Zwiebach}. Furthermore, his idea carries over to homotopy versions of
coalgebras \cite{Saneblidze}, \cite{Smith}, \cite{Umble} and Lie algebras 
\cite{Lada}, and one can deform a classical DG algebra, coalgebra or Lie
algebra to the corresponding homotopy version in a standard way.

This paper introduces the notion of an $A_{\infty }$-bialgebra, which is a
DGM $H$ equipped with \textquotedblleft structurally
compatible\textquotedblright\ operations $\left\{ \omega ^{j,i}:H^{\otimes
i}\rightarrow H^{\otimes j}\right\} _{i,j\geq 1}$ such that $\left( H,\omega
^{1,i}\right) _{i\geq 1}$ is an $A_{\infty }$-algebra and $\left( H,\omega
^{j,1}\right) _{j\geq 1}$ is an $A_{\infty }$-coalgebra. The main result of
this project, the proof of which appears in the sequel \cite{SU3}, is the
fact that over a field, the homology of every $A_{\infty }$-bialgebra
inherits an $A_{\infty }$-bialgebra structure. In particular, the Hopf
algebra structure on a classical Hopf algebra extends to an $A_{\infty }$%
-bialgebra structure and the $A_{\infty }$-bialgebra structure on the
homology of a loop space specializes to the $A_{\infty }$-(co)algebra
structures observed by Gugenheim \cite{Gugenheim} and Kadeishvili \cite%
{Kadeishvili1}. Thus loop space homology provides a primary family of
examples. In fact, one can introduce an $A_{\infty }$-bialgebra structure on
the double cobar construction of H.-J. Baues \cite{Baues1}.

The problem that motivated this project was to classify rational loop spaces
that share a fixed Pontryagin algebra. This problem was considered by the
second author in the mid 1990's as a deformation problem in some large (but
unknown) rational category containing DG\ Hopf algebras. And it was
immediately clear that if such a category exists, it contains objects with
rich higher order structure that specializes to simultaneous $A_{\infty }$%
-algebra and $A_{\infty }$-coalgebra structures. Evidence of this was
presented by the second author at Jim Stasheff's schriftfest (June 1996) in
a talk entitled, \textquotedblleft In search of higher homotopy Hopf
algebras\textquotedblright\ \cite{Umble2}. Given the perspective of this
project, we conjecture that there exists a deformation theory for $A_{\infty
}$-bialgebras in which the infinitesimal deformations of classical DG\
bialgebra's observed in that talk approximate $A_{\infty }$-bialgebras to
first order. Shortly thereafter, the first author used perturbation methods
to solve this classification problem \cite{Saneblidze}. The fact that $%
A_{\infty }$-bialgebras appear implicitly in this solution led to the
collaboration in this project.

Given a DGM $H,$ let $U=End(TH)$ be the associated universal PROP. We
construct internal and external cup products on $C^{\ast}\left( P;\mathbf{U}%
\right) ,$ the cellular chains of permutahedra $P=\sqcup_{n\geq 1}P_{n}$
with coefficients in a certain submodule $\mathbf{U\subset}TTU.$ The first
is defined for every polytope and in particular for each $P_{n};$ the second
is defined globally on $C^{\ast}\left( P;\mathbf{U}\right) $ and depends
heavily on the representation of faces of permutahedra as leveled trees (see
our prequel \cite{SU2}, for example). These cup products give rise to a
biderivative operator $Bd$ on $\mathbf{U}$ with the following property:
Given $\omega\in\mathbf{U},$ there is a unique element $d_{\omega}\in\mathbf{%
U}$ fixed by the action of $Bd$ that bimultiplicatively extends $\omega.$ We
define a (non-bilinear) operation $\circledcirc$ on $\mathbf{U} $ in terms
of $Bd$ and use it to define the notion of an $A_{\infty}$-bialgebra. The
paper is organized as follows:\ Cup products are constructed in Section 2,
the biderivative is defined in Section 3 and $A_{\infty}$-bialgebras are
defined in Section 4.

\section{Cochain Algebras Over the Universal PROP}

Let $R$ be a commutative ring with identity and let $H$ be an $R$-free DGM
of finite type. For $x,y\in \mathbb{N}$, let $U_{y,x}={Hom}\left( H^{\otimes
x},H^{\otimes y}\right) $ and view $U_{H}=End(TH)$ as the bigraded module%
\begin{equation*}
U_{\ast ,\ast }={\bigoplus_{x,y\in \mathbb{N}}U_{y,x}}.
\end{equation*}%
Given matrices $X=\left[ x_{ij}\right] $ and $Y=\left[ y_{ij}\right] \in 
\mathbb{N}^{q\times p},$ consider the module%
\begin{eqnarray*}
U_{Y,X} &=&(U_{y_{11},x_{11}}\otimes \cdots \otimes
U_{y_{1p},x_{1p}})\otimes \cdots \otimes (U_{y_{q1},x_{q1}}\otimes \cdots
\otimes U_{y_{qp},x_{qp}}) \\
&\subset &\left( U^{\otimes p}\right) ^{\otimes q}\subset TTU.
\end{eqnarray*}%
Represent a monomial $A\in U_{Y,X}$ as the $q\times p$ matrix $\left[ A%
\right] =\left[ \theta _{y_{ij},x_{ij}}\right] $ with rows thought of as
elements of $U^{\otimes p}\subset TU$. We refer to $A$ as a $q\times p$
monomial; we often abuse notation and write $A$ when we mean $\left[ A\right]
.$ Note that 
\begin{equation*}
\bigoplus_{_{X,Y\in \mathbb{N}^{q\times p}}}U_{Y,X}=(U^{\otimes p})^{\otimes
q};
\end{equation*}%
in Subsection \ref{subsection1} below we construct the \textquotedblleft
upsilon product\textquotedblright\ on the module 
\begin{equation*}
M=\bigoplus_{_{\substack{ X,Y\in \mathbb{N}^{q\times p}  \\ p,q\geq 1}}%
}U_{Y,X}=\bigoplus_{_{p,q\geq 1}}(U^{\otimes p})^{\otimes q}.
\end{equation*}

In particular, given $\mathbf{x}=\left( x_{1},\ldots ,x_{p}\right) \in 
\mathbb{N}^{p}$ and $\mathbf{y}=\left( y_{1},\ldots ,y_{q}\right) \in 
\mathbb{N}^{q},$ set $X=\left( x_{ij}=x_{j}\right) _{1\leq i\leq q},$ $%
Y=\left( y_{ij}=y_{i}\right) _{1\leq j\leq p}$ and $\mathbf{U}_{\mathbf{x}}^{%
\mathbf{y}}=U_{Y,X}.$ A monomial $A\in \mathbf{U}_{\mathbf{x}}^{\mathbf{y}}$
is represented by a $q\times p$ matrix%
\begin{equation*}
A=\left[ 
\begin{array}{lll}
\theta _{y_{1},x_{1}} & \cdots & \theta _{y_{1},x_{p}} \\ 
\multicolumn{1}{c}{\vdots} & \multicolumn{1}{c}{} & \multicolumn{1}{c}{\vdots
} \\ 
\theta _{y_{q},x_{1}} & \cdots & \theta _{y_{q},x_{p}}%
\end{array}%
\right] .
\end{equation*}%
We refer to the vectors $\mathbf{x}$ and $\mathbf{y}$ as the \emph{%
coderivation }and \emph{derivation} \emph{leaf sequences }of $A\mathbf{,}$
respectively (see Subsection \ref{subsection4}). Note that for $a,b\in 
\mathbb{N}$, monomials in $\mathbf{U}_{a}^{\mathbf{y}}$ and $\mathbf{U}_{%
\mathbf{x}}^{b}$ appear as $q\times 1$ and $1\times p$ matrices. Let%
\begin{equation*}
\mathbf{U}=\bigoplus_{\substack{ \mathbf{x\times y}\in \mathbb{N}^{p}\times 
\mathbb{N}^{q}  \\ p,q\geq 1}}\mathbf{U}_{\mathbf{x}}^{\mathbf{y}}.
\end{equation*}

We graphically represent a monomial $A=\left[ \theta _{y_{j},x_{i}}\right]
\in \mathbf{U}_{\mathbf{x}}^{\mathbf{y}}$ two ways. First as a matrix of
\textquotedblleft double corollas\textquotedblright\ in which entry $\theta
_{y_{j},x_{i}}$ is pictured as two corollas joined at the root--one opening
downward with $x_{i}$ inputs and one opening upward with $y_{j}$
outputs--and second as an arrow in the positive integer lattice $\mathbb{N}%
^{2}$ (see Figure 1). The arrow representation is motivated by the fact that 
$A$ can be thought of as an operator on $\mathbb{N}^{2}.$ Since $H$ has
finite type, $A$ admits a representation as a map%
\begin{equation*}
A:\left( H^{\otimes x_{1}}\otimes \cdots \otimes H^{\otimes x_{p}}\right)
^{\otimes q}\rightarrow \left( H^{\otimes y_{1}}\right) ^{\otimes p}\otimes
\cdots \otimes \left( H^{\otimes y_{q}}\right) ^{\otimes p}.
\end{equation*}%
For $\mathbf{u=}\left( u_{1},\ldots ,u_{k}\right) \in \mathbb{N}^{k},$ let $%
\left\vert \mathbf{u}\right\vert =u_{1}+\cdots +u_{k}$ and identify $\left(
s,t\right) \in \mathbb{N}^{2}$ with the module $\left( H^{\otimes s}\right)
^{\otimes t}.$ Let $\sigma _{s,t}:\left( H^{\otimes s}\right) ^{\otimes t}%
\overset{\approx }{\rightarrow }\left( H^{\otimes t}\right) ^{\otimes s}$ be
the canonical permutation of tensor factors and identify a $q\times p$
monomial $A\in \mathbf{U}_{\mathbf{x}}^{\mathbf{y}}$ with the operator $%
\left( \sigma _{y_{1},p}\otimes \cdots \otimes \sigma _{y_{q},p}\right)
\circ A$ on $\mathbb{N}^{2},$ i.e., the composition%
\begin{equation*}
\left( H^{\otimes \left\vert \mathbf{x}\right\vert }\right) ^{\otimes
q}\approx \left( H^{\otimes x_{1}}\otimes \cdots \otimes H^{\otimes
x_{p}}\right) ^{\otimes q}\overset{A}{\rightarrow }\left( H^{\otimes
y_{1}}\right) ^{\otimes p}\otimes \cdots \otimes \left( H^{\otimes
y_{q}}\right) ^{\otimes p}
\end{equation*}%
\begin{equation*}
\overset{\sigma _{y_{1},p}\otimes \cdots \otimes \sigma _{y_{q},p}}{%
\longrightarrow }\left( H^{\otimes p}\right) ^{\otimes y_{1}}\otimes \cdots
\otimes \left( H^{\otimes p}\right) ^{\otimes y_{q}}\approx \left(
H^{\otimes p}\right) ^{\otimes \left\vert \mathbf{y}\right\vert },
\end{equation*}%
where $\approx $ denotes the canonical isomorphism that changes filtration.
Thus we represent $A$ as an arrow from $\left( \left\vert \mathbf{x}%
\right\vert ,q\right) $ to $\left( p,\left\vert \mathbf{y}\right\vert
\right) $. In particular, a monomial $A\in \mathbf{U}_{a}^{b}$
\textquotedblleft transgresses\textquotedblright\ from $\left( a,1\right) $
to $\left( 1,b\right) .$\vspace*{0.4in}

\hspace*{1.8in}\setlength{\unitlength}{0.001in}%
\begin{picture}
(152,-130)(202,-90) \thicklines

%%%%%%%%%%%%%%%%%%%%%%%%%%%%%%%%%%%%%%%%

%   Corolla_{2,1}

%   Leaves

\put(50,-150){\line( -1,1){150}}
\put(50,-150){\line( 1,1){150}}
\put(50,-150){\line( 0,-1){265}}

%   Vertices

\put(50,-150){\makebox(0,0){$\bullet$}}

%%%%%%%%%%%%%%%%%%%%%%%%%%%%%%%%%%%%%%

%   Corolla_{4,1}

%   Leaves

\put(50,-750){\line( -1,1){150}}
\put(50,-750){\line( 1,1){150}}
\put(50,-750){\line( 1,3){50}}
\put(50,-750){\line( -1,3){50}}
\put(50,-750){\line( 0,-1){265}}

%   Vertices

\put(50,-750){\makebox(0,0){$\bullet$}}

%%%%%%%%%%%%%%%%%%%%%%%%%%%%%%%%%%%%%%

%   Corolla_{2,3}

%   Leaves

\put(550,-150){\line( -1,1){150}}
\put(550,-150){\line( 1,1){150}}
\put(550,-250){\line( -1,-1){150}}
\put(550,-250){\line( 1,-1){150}}
\put(550,-400){\line( 0,1){265}}

%   Vertices

\put(550,-150){\makebox(0,0){$\bullet$}}
\put(550,-250){\makebox(0,0){$\bullet$}}

%%%%%%%%%%%%%%%%%%%%%%%%%%%%%%%%%%%%%%

%   Corolla_{4,3}

%   Leaves

\put(550,-750){\line( -1,1){150}}
\put(550,-750){\line( 1,1){150}}
\put(550,-850){\line( -1,-1){150}}
\put(550,-850){\line( 1,-1){150}}
\put(550,-750){\line( 0,-1){265}}
\put(550,-750){\line( 1,3){50}}
\put(550,-750){\line( -1,3){50}}

%   Vertices

\put(550,-750){\makebox(0,0){$\bullet$}}
\put(550,-850){\makebox(0,0){$\bullet$}}

%%%%%%%%%%%%%%%%%%%%%%%%%%%%%%%%%%%%%%

%   Left bracket

\put(-300,150){\line( 0,-1){1300}}
\put(-300,150){\line( 1,0){100}}
\put(-300,-1150){\line( 1,0){100}}

%  Right bracket

\put(900,150){\line( 0,-1){1300}}
\put(900,150){\line( -1,0){100}}
\put(900,-1150){\line( -1,0){100}}

%%%%%%%%%%%%%%%%%%%%%%%%%%%%%%%%%%%%%%

%   Tag

\put(-900,-500){\makebox(0,0){$A \in \mathbf U_{1,3}^{2,4} \hspace*{0.1in}  \leftrightarrow  \hspace*{0.1in} $}}

\put(1200,-500){\makebox(0,0){$  \leftrightarrow  \hspace*{0.1in} $}}

%%%%%%%%%%%%%%%%%%%%%%%%%%%%%%%%%%%%%

% Arrow diagram

\put(1600,-1150){\line(1,0){1100}}
\put(1600,-1150){\line(0,1){1300}}
\put(1800,-1183){\line(0,1){75}}
\put(2200,-1183){\line(0,1){75}}
\put(1563,-950){\line(1,0){75}}
\put(1563,-150){\line(1,0){75}}
\put(1800,-150){\makebox(0,0){$\bullet$}}
\put(2200,-950){\makebox(0,0){$\bullet$}}
\put(2150,-830){\vector(-1,2){300}}
\put(1450,-950){\makebox(0,0){$2$}}
\put(1450,-150){\makebox(0,0){$6$}}
\put(1800,-1300){\makebox(0,0){$2$}}
\put(2200,-1300){\makebox(0,0){$4$}}

\end{picture}\vspace*{1.4in}

\begin{center}
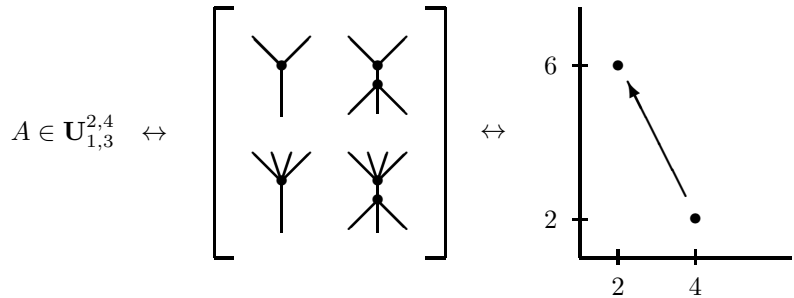
Figure 1. Graphical representations of a typical monomial.\vspace*{0.1in}
\end{center}

\subsection{Products on $\mathbf{U}$\textbf{\label{subsection1}}}

We begin by defining dual associative cross products on $\mathbf{U.}$ Given
a pair of monomials $A\otimes B\in \mathbf{U}_{\mathbf{v}}^{\mathbf{y}%
}\otimes \mathbf{U}_{\mathbf{x}}^{\mathbf{u}},$ define the \emph{wedge} and 
\emph{cech cross products }by 
\begin{equation*}
A\overset{_{\wedge }}{\times }B=\left\{ 
\begin{array}{ll}
A\otimes B, & \text{if }\mathbf{v}=\mathbf{x,} \\ 
0, & \text{otherwise,}%
\end{array}%
\right. \text{ \ and \ }A\overset{_{\vee }}{\times }B=\left\{ 
\begin{array}{ll}
A\otimes B, & \text{if }\mathbf{u}=\mathbf{y,} \\ 
0, & \text{otherwise.}%
\end{array}%
\right.
\end{equation*}%
Then $\mathbf{U}_{\mathbf{x}}^{\mathbf{y}}\overset{_{\wedge }}{\times }%
\mathbf{U}_{\mathbf{x}}^{\mathbf{u}}\subseteq \mathbf{U}_{\mathbf{x}}^{%
\mathbf{y,u}}$ and\ $\mathbf{U}_{\mathbf{v}}^{\mathbf{y}}\overset{_{\vee }}{%
\times }\mathbf{U}_{\mathbf{x}}^{\mathbf{y}}\subseteq \mathbf{U}_{\mathbf{v,x%
}}^{\mathbf{y}};$ denote $\mathbf{\overset{_{\wedge }}{U}=}\left( \mathbf{U,}%
\overset{_{\wedge }}{\times }\right) $ and $\mathbf{\overset{_{\vee }}{U}=}%
\left( \mathbf{U,}\overset{_{\vee }}{\times }\right) $.

Non-zero cross products create block matrices:%
\begin{equation*}
A\overset{_{\wedge }}{\times }B=\QTATOPD[ ] {A}{B}\text{ \ and \ }A\overset{%
_{\vee }}{\times }B=\left[ A\text{ }B\right] .
\end{equation*}%
In terms of arrows, $A\overset{_{\wedge }}{\times }B\in \mathbf{U}_{\mathbf{x%
}}^{\mathbf{y,u}}$ runs from the vertical $x=\left\vert \mathbf{x}%
\right\vert $ to vertical $x=p$ in $\mathbb{N}^{2}$ and $A\overset{_{\vee }}{%
\times }B\in \mathbf{U}_{\mathbf{v,x}}^{\mathbf{y}}$ runs from horizontal $%
y=q$ to $y=\left\vert \mathbf{y}\right\vert .$ Thus an $n\times 1$ monomial $%
A^{\overset{_{\wedge }}{\times }n}\in \mathbf{U}_{a}^{b\cdots b}$ initiates
at $\left( a,n\right) $ and terminates at $\left( 1,nb\right) ;$ a $1\times
n $ monomial $A^{\overset{_{\vee }}{\times }n}\in \mathbf{U}_{a\cdots a}^{b}$
initiates at $\left( na,1\right) $ and terminates at $\left( n,b\right) .$

We also define a composition product on $\mathbf{U.}$

\begin{definition}
A monomial pair $A^{q\times s}\otimes B^{t\times p}=\left[ \theta _{y_{k\ell
},v_{k\ell }}\right] \otimes \left[ \eta _{u_{ij},x_{ij}}\right] \in M%
\mathbf{\otimes }M$ is a

\begin{enumerate}
\item[\textit{(i)}] \underline{Transverse Pair} (TP) if $s=t=1,$ $u_{1,j}=q$
and $v_{k,1}=p$ for all $j,k,$ i.e., setting $x_{j}=x_{1,j}$ and $%
y_{k}=y_{k,1}$ gives 
\begin{equation*}
A\otimes B=\left[ 
\begin{array}{c}
\theta _{y_{1},p} \\ 
\vdots \\ 
\theta _{y_{q},p}%
\end{array}%
\right] \otimes \left[ 
\begin{array}{lll}
\eta _{q,x_{1}} & \cdots & \eta _{q,x_{p}}%
\end{array}%
\right] \in \mathbf{U}_{p}^{\mathbf{y}}\otimes \mathbf{U}_{\mathbf{x}}^{q}.
\end{equation*}

\item[\textit{(ii)}] \underline{Block Transverse Pair} (BTP) if there exist $%
t\times s$ block decompositions $A=\left[ A_{k^{\prime }\ell }^{\prime }%
\right] $ and $B=\left[ B_{ij^{\prime }}^{\prime }\right] $ such that $%
A_{i\ell }^{\prime }\otimes B_{i\ell }^{\prime }$ is a TP for\textit{\ all }$%
i,\ell $.
\end{enumerate}
\end{definition}

\noindent Note that BTP block decomposition is unique; furthermore, $%
A\otimes B\in \mathbf{U}_{\mathbf{v}}^{\mathbf{y}}\otimes \mathbf{U}_{%
\mathbf{x}}^{\mathbf{u}}$ is a BTP if and only if $\mathbf{y}\in \mathbb{N}%
^{|\mathbf{u}|}$ and $\mathbf{x}\in \mathbb{N}^{|\mathbf{v}|}$ if and only
if the initial point of arrow $A$ and the terminal point of arrow $B$
coincide.

\begin{example}
\label{example1}A pairing of monomials $A^{4\times 2}\otimes B^{2\times
3}\in \mathbf{U}_{2,1}^{1,5,4,3}\otimes \mathbf{U}_{1,2,3}^{3,1}$ is a $%
2\times 2$ BTP per the block decompositions\vspace{0.1in}\newline
\hspace*{0.9in}\setlength{\unitlength}{.06in}\special{em:linewidth 0.4pt}%
\linethickness{0.4pt}%
\begin{picture}(118.66,29.34)
\put(5.23,25.00){\makebox(0,0)[cc]{$\theta_{1,2}$}}
\put(5.23,18.00){\makebox(0,0)[cc]{$\theta_{5,2}$}}
\put(5.23,11.67){\makebox(0,0)[cc]{$\theta_{4,2}$}}
\put(5.23,4.67){\makebox(0,0)[cc]{$\theta_{3,2}$}}
\put(1.66,8.67){\dashbox{0.67}(6.33,20.00)[cc]{}}
\put(1.66,2.34){\dashbox{0.67}(6.33,4.67)[cc]{}}
\put(14.00,25.00){\makebox(0,0)[cc]{$\theta_{1,1}$}}
\put(14.00,18.00){\makebox(0,0)[cc]{$\theta_{5,1}$}}
\put(14.00,11.67){\makebox(0,0)[cc]{$\theta_{4,1}$}}
\put(14.00,4.67){\makebox(0,0)[cc]{$\theta_{3,1}$}}
\put(10.33,8.67){\dashbox{0.67}(6.33,20.00)[cc]{}}
\put(10.33,2.34){\dashbox{0.67}(6.33,4.67)[cc]{}}
\put(1.67,15.50){\oval(3.33,28.00)[l]}
\put(16.80,15.50){\oval(3.33,28.00)[r]} \
\put(22.83,16.77){\makebox(0,0)[cc]{\textit{and}}} \
\put(33.00,20.01){\makebox(0,0)[cc]{$\eta_{3,1}$}}
\put(39.58,20.01){\makebox(0,0)[cc]{$\eta_{3,2}$}}
\put(47.53,20.01){\makebox(0,0)[cc]{$\eta_{3,3}$}}
\put(29.67,17.67){\dashbox{0.67}(12.67,5.67)[cc]{}}
\put(44.33,17.67){\dashbox{0.67}(6.67,5.67)[cc]{}}
\put(33.00,12.34){\makebox(0,0)[cc]{$\eta_{1,1}$}}
\put(39.58,12.34){\makebox(0,0)[cc]{$\eta_{1,2}$}}
\put(47.53,13.34){\makebox(0,0)[cc]{$\eta_{1,3}$}}
\put(29.67,10.01){\dashbox{0.67}(12.67,5.67)[cc]{}}
\put(44.53,10.01){\dashbox{0.67}(6.67,5.67)[cc]{}}
\put(28.90,16.27){\oval(2.00,15.33)[l]}
\put(51.56,16.27){\oval(2.00,15.33)[r]}
\put(56.00,16.00){\makebox(0,0)[cc]{$.$}}
\end{picture}\vspace{0.1in}\newline
As arrows, $A$ initializes at $\left( 6,2\right) $ and terminates at $\left(
3,4\right) ;$ $B$ initializes at $\left( 3,4\right) $ and terminates at $%
\left( 2,13\right) .$
\end{example}

When $\mathbf{x\times y}\in \mathbb{N}^{p}\times \mathbb{N}^{q},$ every pair
of monomials $A\otimes B\in \mathbf{U}_{p}^{\mathbf{y}}\otimes \mathbf{U}_{%
\mathbf{x}}^{q}$ is a TP. Define a mapping 
\begin{equation*}
\gamma :\mathbf{U}_{p}^{\mathbf{y}}\otimes \mathbf{U}_{\mathbf{x}%
}^{q}\rightarrow \mathbf{U}_{|\mathbf{x}|}^{|\mathbf{y}|}
\end{equation*}%
by the composition 
\begin{equation*}
\mathbf{U}_{p}^{\mathbf{y}}\otimes \mathbf{U}_{\mathbf{x}}^{q}\overset{\iota
_{q}\otimes \iota _{p}}{\longrightarrow }\mathbf{U}_{pq}^{|\mathbf{y}%
|}\otimes \mathbf{U}_{|\mathbf{x}|}^{qp}\overset{1\otimes \sigma
_{q,p}^{\ast }}{\longrightarrow }\mathbf{U}_{pq}^{|\mathbf{y}|}\otimes 
\mathbf{U}_{|\mathbf{x}|}^{pq}\overset{\circ }{\rightarrow }\mathbf{U}_{|%
\mathbf{x}|}^{|\mathbf{y}|},
\end{equation*}%
where $\iota _{p}$ and $\iota _{q}$ are the canonical isomorphisms. Then for 
$A=\left[ \theta _{y_{k},p}\right] \in \mathbf{U}_{p}^{\mathbf{y}}$ and $B=%
\left[ \eta _{q,x_{j}}\right] \in \mathbf{U}_{\mathbf{x}}^{q},$ we have%
\begin{equation*}
\gamma \left( A\otimes B\right) =\left( \theta _{y_{1},p}\otimes \cdots
\otimes \theta _{y_{q},p}\right) \sigma _{q,p}\left( \eta _{q,x_{1}}\otimes
\cdots \otimes \eta _{q,x_{p}}\right) ;
\end{equation*}%
denote this expression either by $A\cdot B$ or $\gamma (\theta
_{y_{1},p},\ldots ,\theta _{y_{q},p};\eta _{q,x_{1}},\ldots ,\eta
_{q,x_{p}}).$ The $\gamma $-product on matrices of double corollas is
typically a matrix of non-planar graphs (see Figure 2). Note that $\gamma $
agrees with the composition product on the universal preCROC \cite{borya}.

More generally, if $A\otimes B$ is a BTP with block decompositions $A=\left[
A_{i\ell }^{\prime }\right] $ and $B=\left[ B_{i\ell }^{\prime }\right] ,$
define $\gamma \left( A\otimes B\right) _{i\ell }=\gamma \left( A_{i\ell
}^{\prime }\otimes B_{i\ell }^{\prime }\right) .$ Then $\gamma $ sends $%
A^{q\times s}\otimes B^{t\times p}\in \mathbf{U}_{\mathbf{v}}^{\mathbf{y}%
}\otimes \mathbf{U}_{\mathbf{x}}^{\mathbf{u}}$ to a $t\times s$ monomial in $%
\mathbf{U}_{\mathbf{x}^{\prime }}^{\mathbf{y}^{\prime }},$ where $\mathbf{x}%
^{\prime }$ and $\mathbf{y}^{\prime }$ are obtained from $\mathbf{x}$ and $%
\mathbf{y}$ by summing $s$ and $t$ successive coordinate substrings: The
length of the $i^{th}$ substring of $\mathbf{x}$ is the length of the row
matrices in the $i^{th}$ column of $B^{\prime };$ the length of the $\ell
^{th}$ substring of $\mathbf{y}$ is the length of the column matrices in the 
$\ell ^{th}$ row of $A^{\prime }.$ In any case, $\gamma \left( A\otimes
B\right) $ is expressed as an arrow from the initial point of $B$ to the
terminal point of $A.$ \vspace{0.6in}

\hspace*{0.8in}\setlength{\unitlength}{0.0015in}%
\begin{picture}
(152,-130)(202,-90) \thicklines

%Left bracket

\put(-280,50){\line ( 0,-1){1200}}
\put(-280,50){\line ( 1,0){100}}
\put(-280,-1150){\line (1,0){100}}

%%%%%%%%%%%%%%%%%%%%%%%%%%%%%%%%%%%%%%

%   Corolla_{2,1}

%   Leaves

\put(0,-150){\line ( 1,-1){100}}
\put(0,-150){\line ( -1,-1){100}}
\put(0,-150){\line( 0,1){150}}

%   Vertices
\put(0,-150){\makebox(0,0){$\bullet$}}

%%%%%%%%%%%%%%%%%%%%%%%%%%%%%%%%%%%%%%

%   Corolla_{2,1}

%   Leaves

\put(0,-550){\line ( 1,-1){100}}
\put(0,-550){\line ( -1,-1){100}}
\put(0,-550){\line( 0,1){150}}

%   Vertices
\put(0,-550){\makebox(0,0){$\bullet$}}

%%%%%%%%%%%%%%%%%%%%%%%%%%%%%%%%%%%%%%

%   Corolla_{2,1}

%   Leaves
\put(0,-950){\line ( 1,-1){100}}
\put(0,-950){\line ( -1,-1){100}}
\put(0,-950){\line( 0,1){150}}

%   Vertices
\put(0,-950){\makebox(0,0){$\bullet$}}

%%%%%%%%%%%%%%%%%%%%%%%%%%%%%%%%%%%%%%

%Right bracket

\put(280,50){\line ( 0,-1){1200}}
\put(280,50){\line ( -1,0){100}}
\put(280,-1150){\line (-1,0){100}}

%%%%%%%%%%%%%%%%%%%%%%%%%%%%%%%%%%%%%%

\put(500,-550){\makebox(0,0){$\bullet$}}

%%%%%%%%%%%%%%%%%%%%%%%%%%%%%%%%%%%%%%

%Left bracket

\put(700,-300){\line ( 0,-1){500}}
\put(700,-300){\line ( 1,0){100}}
\put(700,-800){\line (1,0){100}}

%%%%%%%%%%%%%%%%%%%%%%%%%%%%%%%%%%%%%%

%   Corolla_{1,3}

%   Leaves
\put(900,-550){\line ( 1,1){100}}
\put(900,-550){\line ( -1,1){100}}
\put(900,-700){\line( 0,1){254}}

%   Vertices
\put(900,-550){\makebox(0,0){$\bullet$}}

%%%%%%%%%%%%%%%%%%%%%%%%%%%%%%%%%%%%%%

%   Corolla_{1,3}

%   Leaves
\put(1200,-550){\line ( 1,1){100}}
\put(1200,-550){\line ( -1,1){100}}
\put(1200,-700){\line( 0,1){254}}

%   Vertices
\put(1200,-550){\makebox(0,0){$\bullet$}}

%%%%%%%%%%%%%%%%%%%%%%%%%%%%%%%%%%%%%%

%Right bracket

\put(1400,-300){\line ( 0,-1){500}}
\put(1400,-300){\line ( -1,0){100}}
\put(1400,-800){\line (-1,0){100}}

%%%%%%%%%%%%%%%%%%%%%%%%%%%%%%%%%%%%%%

\put(1600,-550){\makebox(0,0){$=$}}

%%%%%%%%%%%%%%%%%%%%%%%%%%%%%%%%%%%%%%

%   The graph

\put(1900,-250){\line ( 1,-1){345}}
\put(2300,-650){\line ( 1,-1){270}}
\put(1900,-250){\line ( -1,-1){100}}
\put(1800,-350){\line (0,-1){399}}
\put(1900,-250){\line( 0,1){150}}
\put(2270,-250){\line ( 1,-1){297}}
\put(2270,-250){\line ( -1,-1){165}}
\put(1970,-550){\line (1,1){93}}
\put(1970,-550){\line (0,-1){550}}
\put(2270,-250){\line( 0,1){150}}
\put(2640,-250){\line ( 1,-1){100}}
\put(2740,-350){\line (0,-1){399}}
\put(2640,-250){\line ( -1,-1){165}}
\put(2430,-460){\line ( -1,-1){460}}
\put(2640,-250){\line( 0,1){150}}
\put(1970,-920){\line(-1,1){170}}
\put(2568,-548){\line (0,-1){552}}
\put(2568,-920){\line(1,1){171}}

%   Vertices

\put(2270,-250){\makebox(0,0){$\bullet$}}
\put(1900,-250){\makebox(0,0){$\bullet$}}
\put(2640,-250){\makebox(0,0){$\bullet$}}
\put(1970,-920){\makebox(0,0){$\bullet$}}
\put(2568,-920){\makebox(0,0){$\bullet$}}

\end{picture}\vspace*{1.7in}

\begin{center}
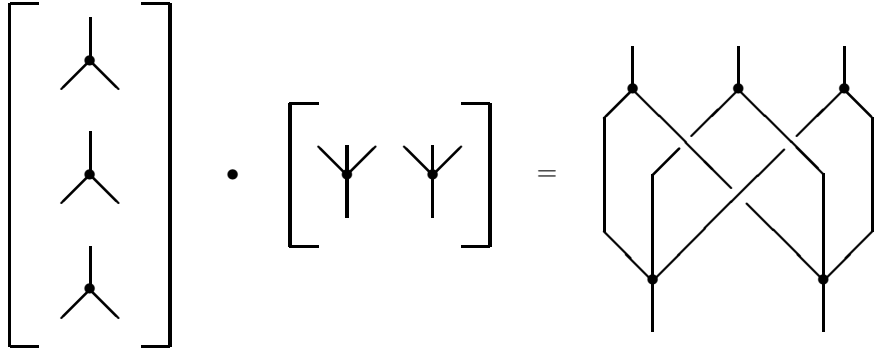
Figure 2. The $\gamma $-product as a non-planar graph.\vspace{0.2in}
\end{center}

Define the \emph{upsilon product }$\Upsilon :M\otimes M\rightarrow M$ on
matrices $A,B\in M$ by 
\begin{equation*}
\Upsilon (A\otimes B)_{i\ell }=\left\{ 
\begin{array}{ll}
\gamma \left( A_{i\ell }^{\prime }\otimes B_{i\ell }^{\prime }\right) , & 
\text{if}\ A\otimes B\ \text{is a}\ \text{BTP} \\ 
&  \\ 
0, & \text{otherwise}%
\end{array}%
\right.
\end{equation*}%
and let $A\cdot B=\Upsilon (A\otimes B).$ Note that $\Upsilon $ restricts to
an associative product on $\mathbf{U}$.

\begin{example}
The $\gamma $-product of the $2\times 2$ BTP $A^{4\times 2}\otimes
B^{2\times 3}$ in Example \ref{example1} is the following $2\times 2$
monomial in $\mathbf{U}_{3,3}^{10,3}:$\vspace{0.2in}\newline
\hspace*{0.1in}\setlength{\unitlength}{1.1mm}\special{em:linewidth 0.4pt}%
\linethickness{0.4pt}%
\begin{picture}(118.66,29.34)
\put(5.20,25.00){\makebox(0,0)[cc]{$\theta_{1,2}$}}
\put(5.20,18.00){\makebox(0,0)[cc]{$\theta_{5,2}$}}
\put(5.20,11.67){\makebox(0,0)[cc]{$\theta_{4,2}$}}
\put(5.20,4.67){\makebox(0,0)[cc]{$\theta_{3,2}$}}
\put(1.66,8.67){\dashbox{0.67}(6.33,20.00)[cc]{}}
\put(1.66,2.34){\dashbox{0.67}(6.33,4.67)[cc]{}}
\put(13.85,25.00){\makebox(0,0)[cc]{$\theta_{1,1}$}}
\put(13.85,18.00){\makebox(0,0)[cc]{$\theta_{5,1}$}}
\put(13.85,11.67){\makebox(0,0)[cc]{$\theta_{4,1}$}}
\put(13.85,4.67){\makebox(0,0)[cc]{$\theta_{3,1}$}}
\put(10.33,8.67){\dashbox{0.67}(6.33,20.00)[cc]{}}
\put(10.33,2.34){\dashbox{0.67}(6.33,4.67)[cc]{}}
\put(1.67,15.34){\oval(3.33,28.00)[l]}
\put(16.80,15.00){\oval(3.33,28.00)[r]} \

\put(20.83,15.77){\makebox(0,0)[cc]{$\cdot$}} \
\put(28.00,20.01){\makebox(0,0)[cc]{$\eta_{3,1}$}}
\put(34.58,20.01){\makebox(0,0)[cc]{$\eta_{3,2}$}}
\put(42.53,20.01){\makebox(0,0)[cc]{$\eta_{3,3}$}}
\put(24.67,17.67){\dashbox{0.67}(12.67,5.67)[cc]{}}
\put(39.33,17.67){\dashbox{0.67}(6.67,5.67)[cc]{}}
\put(28.00,12.34){\makebox(0,0)[cc]{$\eta_{1,1}$}}
\put(34.58,12.34){\makebox(0,0)[cc]{$\eta_{1,2}$}}
\put(42.53,12.34){\makebox(0,0)[cc]{$\eta_{1,3}$}}
\put(24.67,10.01){\dashbox{0.67}(12.67,5.67)[cc]{}}
\put(39.53,10.01){\dashbox{0.67}(6.67,5.67)[cc]{}}
\put(23.90,16.60){\oval(2.00,15.33)[l]}
\put(46.56,16.60){\oval(2.00,15.33)[r]}

\put(51.00,16.00){\makebox(0,0)[cc]{$=$}} \ \
\put(69.32,20.67){\makebox(0,0)[cc]{$_{\gamma(\theta_{1,2},\theta_{5,2},\theta_{4,2};\,\eta_{3,1},\eta_{3,2})}$}}
\put(97.00,20.67){\makebox(0,0)[cc]{$_{\gamma(\theta_{1,1},\theta_{5,1},\theta_{4,1};\,\eta_{3,3})}$}}
\put(69.31,12.67){\makebox(0,0)[cc]{$_{\gamma(\theta_{3,2};\,\eta_{1,1},\eta_{1,2})}$}}
\put(97.00,12.67){\makebox(0,0)[cc]{$_{\gamma(\theta_{3,1};\,\eta_{1,3})}$}}
\put(55.16,16.60){\oval(1.67,15.33)[l]}
\put(108.33,16.60){\oval(2.67,15.33)[r]}
\put(112.00,16.67){\makebox(0,0)[cc]{$.$}}
\end{picture}\vspace{0.1in}\newline
The row matrices in successive columns of the block decomposition of $B$
have respective lengths $2$ and $1;$ thus $\mathbf{x}^{\prime }=\left(
3,3\right) $ is obtained from $\mathbf{x=}\left( 1,2,3\right) .$ Similarly,
the column matrices in successive rows of the block decomposition of $A$
have respective lengths $3$ and $1;$ thus $\mathbf{y}^{\prime }=\left(
10,3\right) $ is obtained from $\mathbf{y}=\left( 1,5,4,3\right) .$ Finally,
the map $A\cdot B:\left( H^{\otimes 3}\otimes H^{\otimes 3}\right) ^{\otimes
2}\rightarrow \left( H^{\otimes 10}\right) ^{\otimes 2}\otimes \left(
H^{\otimes 3}\right) ^{\otimes 2}$ is expressed as an arrow initializing at $%
\left( 6,2\right) $ and terminating at $\left( 2,13\right) .$
\end{example}

\subsection{Cup products on $C^{\ast }\left( P,\mathbf{U}\right) $}

Let $C_{\ast }\left( X\right) $ denote the cellular chains on a polytope $X$
and assume that $C_{\ast }\left( X\right) $ comes equipped with a diagonal $%
\Delta _{X}:C_{\ast }\left( X\right) \rightarrow C_{\ast }\left( X\right)
\otimes C_{\ast }\left( X\right) .$ Let $G$ be a module (graded or
ungraded); if $G$ is graded, ignore the grading and view $G$ as a graded
module concentrated in degree zero. The cellular $k$-cochains on $X$\ with
coefficients in $G$ is the graded module%
\begin{equation*}
C^{k}\left( X;G\right) =Hom^{-k}\left( C_{\ast }\left( X\right) ,G\right) .
\end{equation*}%
When $G$ is a DGA with multiplication $\mu ,$ the diagonal $\Delta _{X}$
induces a DGA\ structure on $C^{\ast }\left( X;G\right) $ with cup product 
\begin{equation*}
f\smile g=\mu \left( f\otimes g\right) \Delta _{X}.
\end{equation*}%
Unless explicitly indicated otherwise, non-associative cup products with
multiple factors are parenthesized on the extreme left, i.e., $f\smile
g\smile h=\left( f\smile g\right) \smile h.$

In our prequel \cite{SU2} we constructed an explicit non-coassociative
non-cocommuta- tive diagonal $\Delta _{P}$ on the cellular chains of
permutahedra $C_{\ast }\left( P_{n}\right) $ for each $n\geq 1.$ Thus we
immediately obtain non-associative, non-commutative DGA's $C^{\ast }(P_{n};%
\overset{_{\wedge }}{\mathbf{U}})$ and $C^{\ast }(P_{n};\overset{_{\vee }}{%
\mathbf{U}})$ with respective wedge and cech cup products $\wedge $ and $%
\vee $. Of course, summing over all $n$ gives wedge and cech cup products on 
$C^{\ast }(P;\mathbf{U}).$

The modules $C^{\ast }(P;\overset{_{\wedge }}{\mathbf{U}})$ and $C^{\ast }(P;%
\overset{_{\vee }}{\mathbf{U}})$ are equipped with second cup products $%
\wedge _{\ell }$ and $\vee _{\ell },$ which arise from the $\Upsilon $%
-product on $\mathbf{U}$ together with the \textquotedblleft level
coproduct.\textquotedblright\ Recall that $m$-faces of $P_{n+1}$ are indexed
by PLT's with $n+2$ leaves, $n-m+1$ levels and root in level $n-m+1$ (see 
\cite{Loday} or \cite{SU2}, for example). The \emph{level coproduct} 
\begin{equation*}
\Delta _{\ell }:C_{\ast }(P)\rightarrow C_{\ast }(P)\otimes C_{\ast }(P)
\end{equation*}%
vanishes on $e^{n}\subset P_{n+1}$ and is defined on proper $m$-faces $e^{m}$
as follows: For each $k,$ prune the tree of $e^{m}$ between levels $k$ and $%
k+1$ and sequentially number the stalks or trees removed from left-to-right.
Let $e_{k}^{\prime }$ denote the pruned tree; let $e_{k}^{\prime \prime }$
denote the tree obtained by attaching all stalks and trees removed during
pruning to a common root (see Figure 4). Then 
\begin{equation*}
\Delta _{\ell }\left( e^{m}\right) =\sum\limits_{1\leq k\leq
n-m}e_{k}^{\prime }\otimes e_{k}^{\prime \prime }.
\end{equation*}%
\vspace{0.3in}

\hspace*{0.2in}\setlength{\unitlength}{0.04in}{%
\begin{picture}(112.34,25.33)
\thicklines
%%%%%% tree #1
%%%%%
\put(8.0,20.33){\line(-1,2){2.40}}
\put(8.0,20.33){\makebox(0,0){$%
\bullet$}}
\put(8.0,20.33){\line(1,2){2.40}}
\put(15,20.33){%
\line(-1,2){2.40}}
\put(15,20.33){\makebox(0,0){$\bullet$%
}}
\put(15,20.33){\line(1,2){2.40}}
\put(5.6,15.33){\line(-1,2){4.95}}
\put(5.6,15.33){\makebox(0,0){$\bullet$}}
\put(5.6,15.33){%
\line(1,2){2.40}}
\put(15.02,20.00){\line(0,-1){5.00}}
\put(15.02,15.33){%
\makebox(0,0){$\bullet$}}
\put(10.31,10.5){\line(-1,1){4.67}}
\put(10.31,10.5){\line(1,1){4.67}}
\put(10.31,10.5){\makebox(0,0){$\bullet$%
}}
\put(10.31,10.5){\line(0,-1){4.67}}
%%%%%% tree #2
%%%%%%
\put(38.5,10.5){\line(-1,1){4.67}}
\put(38.5,10.5){%
\line(1,1){4.67}}
\put(38.5,10.5){\makebox(0,0){$\bullet$%
}}
\put(33.67,15.33){\line(-1,2){2.00}}
\put(33.67,15.33){%
\line(1,2){2.00}}
\put(33.67,15.33){\makebox(0,0){$\bullet$%
}}
\put(43.33,15.33){\makebox(0,0){$\bullet$}}
\put(38.5,10.5){%
\line(0,-1){4.67}}
\put(43.33,15.33){\line(0,1){4.0}}
%%%%%% tree #3
%%%%%%
\put(58.0,20.33){\line(-1,2){2.40}}
\put(58.0,20.33){\makebox(0,0){$%
\bullet$}}
\put(58.0,20.33){\line(1,2){2.40}}
\put(64.5,20.33){%
\line(-1,2){2.4}}
\put(64.5,20.33){\makebox(0,0){$\bullet$%
}}
\put(64.5,20.33){\line(1,2){2.40}}
\put(58.0,6.00){\line(0,1){15.0}}
\put(58.0,10.5){\makebox(0,0){$\bullet$}}
\put(57.8,10.5){\line(-1,2){7.4}}
\put(58.0,10.5){\line(2,3){6.8}}
%%%%% tree #4
%%%%%%%
\put(82.01,10.5){\line(-1,1){4}}
\put(82.01,10.5){\line(1,1){4}}
\put(82.01,10.5){\makebox(0,0){$\bullet$%
}}
\put(82.01,10.5){\line(0,-1){4.67}}
%%%%% tree #5
%%%%%%%
\put(101.03,20.33){\line(-1,2){2.40}}
\put(101.03,20.33){%
\makebox(0,0){$\bullet$}}
\put(101.03,20.33){\line(1,2){2.40}}
\put(108.03,20.33){\line(-1,2){2.40}}
\put(108.03,20.33){\makebox(0,0){$%
\bullet$}}
\put(108.03,20.33){\line(1,2){2.40}}
\put(98.63,15.33){%
\line(-1,2){4.95}}
\put(98.63,15.33){\makebox(0,0){$\bullet$%
}}
\put(98.63,15.33){\line(1,2){2.40}}
\put(108.05,20.00){%
\line(0,-1){5.00}}
\put(108.05,15.33){\makebox(0,0){$\bullet$%
}}
\put(103.34,10.5){\line(-1,1){4.67}}
\put(103.34,10.5){%
\line(1,1){4.67}}
\put(103.34,10.5){\makebox(0,0){$\bullet$%
}}
\put(103.34,10.5){\line(0,-1){4.67}}
%%%%%% labels and symbols
%%%%%%%%
\put(-3.10,20.33){\makebox(0,0)[cc]{$_{k=1}$}}
\put(-3.10,15.33){%
\makebox(0,0)[cc]{$_{k=2}$}}
\put(-3.10,10.50){\makebox(0,0)[cc]{$_{k=3}$%
}}
\put(23.34,15.66){\makebox(0,0)[cc]{$\Delta_{\ell}$}}
\put(50.00,13.66){\makebox(0,0)[cc]{$\otimes$}}
\put(92.01,13.33){%
\makebox(0,0)[cc]{$\otimes$}}
\put(70.34,13.33){\makebox(0,0)[cc]{$+$%
}}
\put(19.00,12.33){\vector(1,0){8.33}}
\end{picture}}

\begin{center}
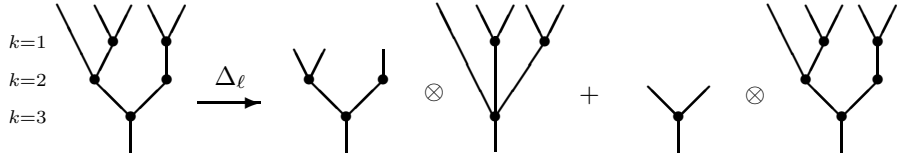
Figure 3: $\Delta _{\ell }(24|1|3)=1|2\otimes 24|13+1\otimes 24|1|3.$\vspace{%
0.15in}
\end{center}

\noindent Obviously, $\Delta _{\ell }$ is non-counital, non-cocommutative
and non-coassociative; in fact, it fails to be a chain map. Fortunately this
is not an obstruction to lifting the $\gamma $-product on $\mathbf{U}$ to a $%
\smile _{\ell }$-product on $C^{\ast }(P;\mathbf{U})$ since we restrict to
certain canonically associative subalgebras of $\mathbf{U}$. For $\varphi
,\varphi ^{\prime }\in C^{\ast }(P;\mathbf{\overset{_{\wedge }}{M}})$ define 
$\varphi \wedge _{\ell }\varphi ^{\prime }=\varphi \smile _{\ell }\varphi
^{\prime }$ and for $\psi ,\psi ^{\prime }\in C^{\ast }(P;\mathbf{\overset{%
_{\vee }}{M}})$ define $\psi \vee _{\ell }\psi ^{\prime }=\psi ^{\prime
}\smile _{\ell }\psi $. Some typical $\wedge _{\ell }$-products appear in
Example \ref{ex1} below.

\subsection{Leaf sequences\label{subsection4}}

Let $T$ be a PLT with at least 2 leaves. Prune $T$ immediately below the
first level, trimming off $k$ stalks and corollas. Number them sequentially
from left-to-right and let $n_{j}$ denote the number of leaves in the $%
j^{th} $ corolla (if $T$ is a corolla, $k=1$ and the pruned tree is a
stalk). The \emph{leaf sequence} of $T$ is the vector $\left( n_{1},\ldots
,n_{k}\right) \in \mathbb{N}^{k}.$

Given integers $n$ and $k$ with $1\leq k\leq n+1,$ let $\mathbf{n=}\left(
n_{1},\ldots ,n_{k}\right) \in \left\{ \mathbf{x}\in \mathbb{N}%
^{k}\right\vert $ $|\mathbf{x}|=n+2\}.\,$When $k=1,$ $e_{\mathbf{n}}$
denotes the $\left( n+2\right) $-leaf corolla. Otherwise, $e_{\mathbf{n}}$
denotes the 2-levelled tree with leaf sequence $\mathbf{n}.$ Now consider
the DGA $\mathbf{U}$ with its $\gamma $-product. Given a codim $0$ or $1$
face $e_{\mathbf{n}}\subset P$ and a cochain $\varphi \in C^{\ast }(P;%
\mathbf{U}),$ let $\varphi _{\mathbf{n}}=\varphi (e_{\mathbf{n}}).$

\begin{example}
\label{ex1}Let $\varphi \in C^{0}(P;\overset{\wedge }{\mathbf{U}})$ and $%
\bar{\varphi}\in C^{1}(P;\overset{\wedge }{\mathbf{U}}).$ When $n=1,$ the
proper faces of $P_{2}$ are its vertices $1|2$ and $2|1$ with $\Delta _{\ell
}(1|2)=1\otimes 1|2$ and $\Delta _{\ell }(2|1)=1\otimes 2|1$. Evaluating $%
\wedge _{\ell }$-squares on vertices gives the compositions 
\begin{equation*}
\begin{array}{l}
\varphi ^{2}(1|2)=\varphi _{2}\varphi _{21}%
\end{array}%
\hspace{0.5in}%
\begin{array}{l}
\varphi ^{2}(2|1)=\varphi _{2}\varphi _{12}.%
\end{array}%
\end{equation*}%
When $n=2,$ the proper faces of $P_{3}$ are its edges and vertices (see
Figure 4). Evaluating quadratic and cubic $\wedge _{\ell }$-products on
edges and vertices gives 
\begin{equation*}
\begin{array}{lll}
{\bar{\varphi}}^{2}(1|23) & = & \bar{\varphi}_{3}\bar{\varphi}_{211} \\ 
\bar{\varphi}^{2}(2|13) & = & \bar{\varphi}_{3}\bar{\varphi}_{121} \\ 
\bar{\varphi}^{2}(3|12) & = & \bar{\varphi}_{3}\bar{\varphi}_{112} \\ 
\varphi \bar{\varphi}(12|3) & = & \varphi _{2}\bar{\varphi}_{31} \\ 
\varphi \bar{\varphi}(13|2) & = & \varphi _{2}\bar{\varphi}_{22} \\ 
\varphi \bar{\varphi}(23|1) & = & \varphi _{2}\bar{\varphi}_{13}%
\end{array}%
\hspace{0.2in}%
\begin{array}{lll}
\varphi ^{2}\bar{\varphi}(1|2|3) & = & \varphi _{2}\varphi _{21}\bar{\varphi}%
_{211} \\ 
\varphi ^{2}\bar{\varphi}(1|3|2) & = & \varphi _{2}\varphi _{12}\bar{\varphi}%
_{211} \\ 
\varphi ^{2}\bar{\varphi}(2|1|3) & = & \varphi _{2}\varphi _{21}\bar{\varphi}%
_{121} \\ 
\varphi ^{2}\bar{\varphi}(2|3|1) & = & \varphi _{2}\varphi _{12}\bar{\varphi}%
_{121} \\ 
\varphi ^{2}\bar{\varphi}(3|1|2) & = & \varphi _{2}\varphi _{21}\bar{\varphi}%
_{112} \\ 
\varphi ^{2}\bar{\varphi}(3|2|1) & = & \varphi _{2}\varphi _{12}\bar{\varphi}%
_{112}.%
\end{array}%
\end{equation*}%
\vspace*{0.1in}
\end{example}

\begin{center}
\setlength{\unitlength}{0.0005in}%
\begin{picture}
(2975,2685)(3126,-2038) \thicklines\put(3601,239){\line( 1,
0){1800}} \put(5401,239){\line( 0,-1){1800}}
\put(5401,-1561){\line(-1, 0){1800}} \put(3601,-1561){\line( 0,
1){1800}} \put(3601,239){\makebox(0,0){$\bullet$}}
\put(3601,-661){\makebox(0,0){$\bullet$}}
\put(3601,-1561){\makebox (0,0){$\bullet$}}
\put(5401,239){\makebox(0,0){$\bullet$}} \put
(5401,-661){\makebox(0,0){$\bullet$}}
\put(5401,-1561){\makebox(0,0){$\bullet $}}
\put(4500,-680){\makebox(0,0){$123$}} \put(2800,-1861){\makebox
(0,0){$1|2|3$}} \put(2800,-699){\makebox(0,0){$ 1|3|2$}}
\put(2800,464){\makebox(0,0){$3|1|2$}}
\put(6200,-1861){\makebox(0,0){$2|1|3$}} \put
(6200,-699){\makebox(0,0){$2|3|1$}} \put(6200,464){\makebox
(0,0){$3|2|1$}} \put(3110,-1260){\makebox(0,0){$1|23$}}
\put(4550,530){\makebox(0,0){$3|12$}} \put(3100,-111){\makebox
(0,0){$13|2$}} \put(5900,-111){\makebox(0,0){$23|1$}}
\put(5900,-1260){\makebox(0,0){$2|13$}} \put(4550,-1890){\makebox
(0,0){$12|3$}}
\end{picture}\vspace{0.1in}

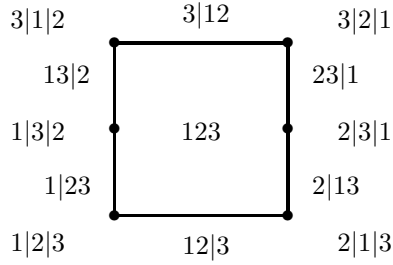
Figure 4: The permutahedron $P_{3}.$
\end{center}

\section{The biderivative}

The definition of the biderivative operator $Bd:\mathbf{U}\rightarrow 
\mathbf{U}$ requires some notational preliminaries. Let $\mathbf{x}%
_{i}\left( r\right) =\left( 1,\ldots ,r,\ldots ,1\right) $ with $r\geq 1$ in
the $i^{th}$ position; the subscript $i$ will be suppressed unless we need
its precise value; in particular, let $\boldsymbol{1}^{k}=\mathbf{x}\left(
1\right) \in \mathbb{N}^{k}.$ Again, we often suppress the superscript and
write $\boldsymbol{1}$ when the context is clear. Let 
\begin{equation*}
\mathbf{U}_{0}={\bigoplus\limits_{\mathbf{x,y=}\boldsymbol{1}}}\mathbf{U}_{%
\mathbf{x}}^{\mathbf{y}}\text{ \ \ and \ \ }\mathbf{U}_{+}=\mathbf{U\diagup U%
}_{0}={\bigoplus\limits_{\mathbf{x\neq }\boldsymbol{1}\text{ or }\mathbf{%
y\neq }\boldsymbol{1}}}\mathbf{U}_{\mathbf{x}}^{\mathbf{y}};
\end{equation*}%
also denote the submodules%
\begin{equation*}
\begin{array}{llllllll}
\mathbf{U}_{u_{0}} & = & {\dbigoplus\limits_{\mathbf{x}\in \mathbb{N}^{p};%
\text{ }|\mathbf{x}|>p\geq 1}}\mathbf{U}_{\mathbf{x}}^{1} &  &  & \mathbf{U}%
_{v_{0}} & = & {\dbigoplus\limits_{\mathbf{y}\in \mathbb{N}^{q};\text{ }|%
\mathbf{y}|>q\geq 1}}\mathbf{U}_{1}^{\mathbf{y}} \\ 
&  &  &  &  &  &  &  \\ 
\mathbf{U}_{u} & = & {\dbigoplus\limits_{\substack{ \mathbf{x}\in \mathbb{N}%
^{p};\text{ }|\mathbf{x}|>1  \\ p,q\geq 1}}}\mathbf{U}_{\mathbf{x}}^{q} &  & 
& \mathbf{U}_{v} & = & {\dbigoplus\limits_{_{\substack{ \mathbf{y}\in 
\mathbb{N}^{q};\text{ }|\mathbf{y}|>1  \\ p,q\geq 1}}}}\mathbf{U}_{p}^{%
\mathbf{y}}%
\end{array}%
\end{equation*}%
and note that%
\begin{equation*}
\mathbf{U}_{u\cap v}=\mathbf{U}_{u}\cap \mathbf{U}_{v}=\bigoplus_{p,q\geq 2}%
\mathbf{U}_{p}^{q}.
\end{equation*}%
Monomials in $\mathbf{U}_{u}$ and $\mathbf{U}_{v}$ are respectively row and
column matrices. In terms of arrows, $\mathbf{U}_{0}$ consists of all arrows
of length zero; $\mathbf{U}_{+}$ consists of all arrows of positive length.
Arrows in $\mathbf{U}_{u}$ initiate on the $x$-axis at $\left( \left\vert 
\mathbf{x}\right\vert ,1\right) ,$ $\left\vert \mathbf{x}\right\vert >1,$
and terminate in the region $x\leq \left\vert \mathbf{x}\right\vert ;$ in
particular, arrows in $\mathbf{U}_{u_{0}}$ lie on the $x$-axis and terminate
at $\left( p,1\right) .$ Arrows in $\mathbf{U}_{v}$ initiate in the region $%
y\leq \left\vert \mathbf{y}\right\vert $ and terminate at $\left(
1,\left\vert \mathbf{y}\right\vert \right) ,$ $\left\vert \mathbf{y}%
\right\vert >1;$ in particular, arrows in $\mathbf{U}_{v_{0}}$ lie on the $y$%
-axis and initiate at $\left( 1,q\right) .$ Thus arrows in $\mathbf{U}%
_{u\cap v}$ \textquotedblleft transgress\textquotedblright\ from the $x$ to
the $y$ axis.

\subsection{The non-linear operator $BD$}

Recall that $\mathbf{n}$ is a leaf sequence if and only if $\mathbf{n\neq 
\boldsymbol{1;}}$ when this occurs, $e_{\mathbf{n}}$ is a face of $%
P_{\left\vert \mathbf{n}\right\vert -1}$ in dimension $\left\vert \mathbf{n}%
\right\vert -2$ or $\left\vert \mathbf{n}\right\vert -3.$ Let 
\begin{equation*}
\overset{_{\wedge }}{T}op(\mathbf{U}_{+})\subset C^{\ast
}(P;\tbigoplus\nolimits_{\mathbf{x\neq \boldsymbol{1}}\text{ or }\mathbf{%
y\neq \boldsymbol{1}}}\mathbf{U}_{\mathbf{x}}^{\mathbf{y}})
\end{equation*}%
be the submodule supported on $e_{\mathbf{x}}$ when $\mathbf{x}\neq 
\boldsymbol{1}$ or on $e_{\mathbf{y}}$ otherwise. Dually, let 
\begin{equation*}
\overset{_{\vee }}{T}op(\mathbf{U}_{+})\subset C^{\ast
}(P;\tbigoplus\nolimits_{\mathbf{x\neq \boldsymbol{1}}\text{ or }\mathbf{%
y\neq \boldsymbol{1}}}\mathbf{U}_{\mathbf{x}}^{\mathbf{y}})
\end{equation*}%
be the submodule supported on $e_{\mathbf{y}}$ when $\mathbf{y}\neq 
\boldsymbol{1}$ or on $e_{\mathbf{x}}$ otherwise. When $\mathbf{x,y}\neq 
\boldsymbol{1,}$ a monomial $A\in \mathbf{U}_{\mathbf{x}}^{\mathbf{y}}$ is
identified with the cochains $\varphi _{A}\in \overset{_{\wedge }}{T}op(%
\mathbf{U}_{+})$ and $\psi _{A}\in \overset{_{\vee }}{T}op(\mathbf{U}_{+})$
respectively supported on the codim $0$ or $1$ faces of $P$ with leaf
sequences $\mathbf{x}$ and $\mathbf{y}$ (see Figure 5). Let 
\begin{equation*}
\overset{_{\wedge }}{\pi }:C^{\ast }(P;{\mathbf{U}}_{+})\rightarrow \overset{%
_{\wedge }}{T}op(\mathbf{U}_{+})\text{ \ and \ }\overset{_{\vee }}{\pi }%
:C^{\ast }(P;{\mathbf{U}}_{+})\rightarrow \overset{_{\vee }}{T}op(\mathbf{U}%
_{+})
\end{equation*}%
be the canonical projections.\vspace{1.6in}

\hspace*{2.8in}\setlength{\unitlength}{0.001in}%
\begin{picture}
(152,-130)(202,-90) \thicklines

%%%%%%%%%%%%%%%%%%%%%%%%%%%%%%%%%%%%%%%%

%   Corolla_{2,1}

%   Leaves

\put(50,-150){\line( -1,1){150}}
\put(50,-150){\line( 1,1){150}}
\put(50,-150){\line( 0,-1){265}}

%   Vertices

\put(50,-150){\makebox(0,0){$\bullet$}}

%%%%%%%%%%%%%%%%%%%%%%%%%%%%%%%%%%%%%%

%   Corolla_{4,1}

%   Leaves

\put(50,-750){\line( -1,1){150}}
\put(50,-750){\line( 1,1){150}}
\put(50,-750){\line( 1,3){50}}
\put(50,-750){\line( -1,3){50}}
\put(50,-750){\line( 0,-1){265}}

%   Vertices

\put(50,-750){\makebox(0,0){$\bullet$}}

%%%%%%%%%%%%%%%%%%%%%%%%%%%%%%%%%%%%%%

%   Corolla_{2,3}

%   Leaves

\put(550,-150){\line( -1,1){150}}
\put(550,-150){\line( 1,1){150}}
\put(550,-250){\line( -1,-1){150}}
\put(550,-250){\line( 1,-1){150}}
\put(550,-400){\line( 0,1){265}}

%   Vertices

\put(550,-150){\makebox(0,0){$\bullet$}}
\put(550,-250){\makebox(0,0){$\bullet$}}

%%%%%%%%%%%%%%%%%%%%%%%%%%%%%%%%%%%%%%

%   Corolla_{4,3}

%   Leaves

\put(550,-750){\line( -1,1){150}}
\put(550,-750){\line( 1,1){150}}
\put(550,-850){\line( -1,-1){150}}
\put(550,-850){\line( 1,-1){150}}
\put(550,-750){\line( 0,-1){265}}
\put(550,-750){\line( 1,3){50}}
\put(550,-750){\line( -1,3){50}}

%   Vertices

\put(550,-750){\makebox(0,0){$\bullet$}}
\put(550,-850){\makebox(0,0){$\bullet$}}

%%%%%%%%%%%%%%%%%%%%%%%%%%%%%%%%%%%%%%

%   Left bracket

\put(-300,150){\line( 0,-1){1300}}
\put(-300,150){\line( 1,0){100}}
\put(-300,-1150){\line( 1,0){100}}

%  Right bracket

\put(900,150){\line( 0,-1){1300}}
\put(900,150){\line( -1,0){100}}
\put(900,-1150){\line( -1,0){100}}

%%%%%%%%%%%%%%%%%%%%%%%%%%%%%%%%%%%%%%

%   Codim 1 face e_{2,4}

%   Leaves

\put(270,1200){\line( -1,1){150}}
\put(70,1000){\line( 1,1){350}}
\put(270,1200){\line( 1,3){50}}
\put(270,1200){\line( -1,3){50}}
\put(70,1000){\line( -1,1){350}}
\put(70,1000){\line( 0,-1){300}}
\put(-130,1200){\line(1,1){150}}

%   Vertices

\put(270,1200){\makebox(0,0){$\bullet$}}
\put(70,1000){\makebox(0,0){$\bullet$}}
\put(-130,1200){\makebox(0,0){$\bullet$}}

%  Arrow

\put(70,500){\vector( 0,-1){200}}
\put(500,1000){\makebox(0,0){$e_{2,4}$}}
\put(270,400){\makebox(0,0){$\psi_A$}}

%%%%%%%%%%%%%%%%%%%%%%%%%%%%%%%%%%%%%%

%   Codim 1 face e_{1,3}

%   Leaves

\put(-1250,-350){\line( -1,-1){150}}
\put(-1450,-150){\line( 1,-1){350}}
\put(-1250,-350){\line( 0,-1){170}}
\put(-1450,-150){\line( -1,-1){350}}
\put(-1450,-150){\line(0,1){350}}

%   Vertices

\put(-1250,-350){\makebox(0,0){$\bullet$}}
\put(-1450,-150){\makebox(0,0){$\bullet$}}

%  Arrow

\put(-700,-150){\vector( 1,0){200}}
\put(-1000,-150){\makebox(0,0){$e_{1,3}$}}
\put(-600,-350){\makebox(0,0){$\varphi_A$}}

%%%%%%%%%%%%%%%%%%%%%%%%%%%%%%%%%%%%%%

% Tag

\put(1200,-500){\makebox(0,0){$  \leftrightarrow  \hspace*{0.1in} A$}}

\end{picture}\vspace{1.2in}

\begin{center}
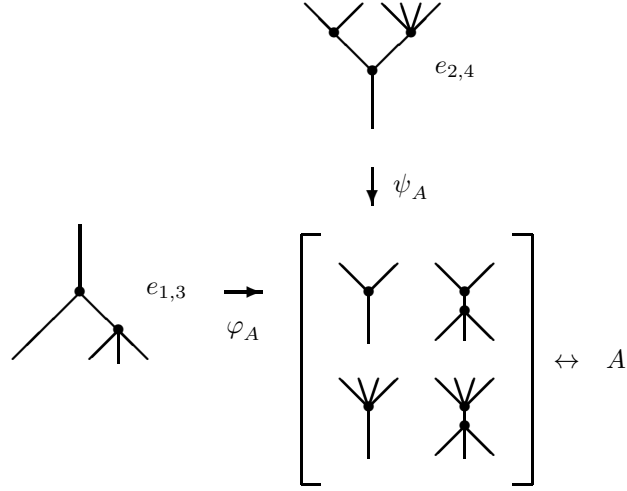
Figure 5: The monomial $A$ is identified with $\varphi _{A}$ and $\psi _{A}$.%
\vspace{0.2in}
\end{center}

For $\mathbf{x}_{i}(n),\mathbf{y}_{i}(n)\in \mathbb{N}^{q},$ let%
\begin{equation*}
\left[ \theta _{1,n}\right] _{i}^{\overset{_{\vee }}{}}=[Id\text{ }\cdots
\theta _{1,n}\cdots \text{ }Id]\in \mathbf{U}_{\mathbf{x}_{i}(n)}^{1}\subset 
\mathbf{U}_{u_{0}}
\end{equation*}%
and%
\begin{equation*}
\left[ \theta _{n,1}\right] _{i}^{\overset{_{\wedge }}{}}=\left[ 
\begin{array}{c}
Id \\ 
\vdots \\ 
\theta _{n,1} \\ 
\vdots \\ 
Id%
\end{array}%
\right] \in \mathbf{U}_{1}^{\mathbf{y}_{i}(n)}\subset \mathbf{U}_{v_{0}}.
\end{equation*}%
Given $\phi \in C^{\ast }(P;\mathbf{U}_{+})$ and $n\geq 2,$ consider the top
dimensional cell $e_{n}\subseteq P_{n-1}$ and components $\phi _{1,n}\left(
e_{n}\right) \in \mathbf{U}_{n}^{1}\subset \mathbf{U}_{u_{0}}$ and $\phi
_{n,1}\left( e_{n}\right) \in \mathbf{U}_{1}^{n}\subset \mathbf{U}_{v_{0}}$
of $\phi \left( e_{n}\right) .$ The \emph{coderivation cochain of }$\phi $
is the global cochain $\phi ^{c}\in \overset{_{\wedge }}{T}op(\mathbf{U}%
_{+}) $ given by%
\begin{equation*}
\phi ^{c}(e_{\mathbf{x}})=\left\{ 
\begin{array}{ll}
\left[ \phi _{1,n}(e_{n})\right] _{i}^{\overset{_{\vee }}{}}, & \text{if}\ 
\mathbf{x}=\mathbf{x}_{i}(n),\text{ }1\leq i\leq q,\text{ }n\geq 2 \\ 
0, & \text{otherwise.}%
\end{array}%
\right.
\end{equation*}
Dually, the \emph{derivation cochain of} $\phi $ is the cochain $\phi
^{a}\in \overset{_{\vee }}{T}op(\mathbf{U}_{+})$ given by%
\begin{equation*}
\phi ^{a}(e_{\mathbf{y}})=\left\{ 
\begin{array}{ll}
\left[ \phi _{n,1}(e_{n})\right] _{i}^{\overset{_{\wedge }}{}}, & \text{if}\ 
\mathbf{y}=\mathbf{y}_{i}(n),\text{ }1\leq i\leq q,\text{ }n\geq 2 \\ 
0, & \text{otherwise.}%
\end{array}%
\right.
\end{equation*}%
Thus $\phi ^{c}$ is supported on the union of the $e_{\mathbf{x}_{i}\left(
n\right) }$'s and takes the value%
\begin{equation*}
\phi ^{c}C_{\ast }\left( P\right) =\dsum_{\substack{ 1\leq i\leq q  \\ q\geq
1 }}[Id\cdots \underset{i^{th}}{\underbrace{\phi _{u_{0}}(e_{n})}}\cdots
Id]^{1\times q}\in \mathbf{U}_{u_{0}},
\end{equation*}%
and dually for $\phi ^{a}$.

Finally, define an operator $\tau :C^{\ast }(P;\mathbf{U_{+}})\rightarrow
C^{\ast }(P;\mathbf{U_{+}})$ on a cochain $\xi \in C^{\ast }(P;\mathbf{U}_{%
\mathbf{x}}^{\mathbf{y}})$ by 
\begin{equation*}
\tau (\xi )(e)=\left\{ 
\begin{array}{ll}
\xi \left( e_{\mathbf{x}}\right) , & \text{if}\ e=e_{\mathbf{y}};\text{ }%
\mathbf{x,y\neq \boldsymbol{1}}, \\ 
\xi \left( e_{\mathbf{y}}\right) , & \text{if}\ e=e_{\mathbf{x}};\text{ }%
\mathbf{x,y\neq \boldsymbol{1}}, \\ 
0, & \text{otherwise}.%
\end{array}%
\right.
\end{equation*}%
Note that $\tau $ is involutory on $\overset{_{\wedge }}{T}op(\mathbf{U}%
_{+})\cap \overset{_{\vee }}{T}op(\mathbf{U}_{+}).$

We are ready to define the non-linear operator $BD.$ First define operators 
\begin{equation*}
\overset{_{\wedge }}{B}D:C^{\ast }(P;\overset{_{\wedge }}{\mathbf{U}}%
_{+})\rightarrow C^{\ast }(P;\overset{_{\wedge }}{\mathbf{U}}_{+})\ \ \text{%
and \ }\overset{_{\vee }}{B}D:C^{\ast }(P;\overset{_{\vee }}{\mathbf{U}}%
_{+})\rightarrow C^{\ast }(P;\overset{_{\vee }}{\mathbf{U}}_{+})
\end{equation*}%
by%
\begin{equation*}
\overset{_{\wedge }}{B}D(\varphi )=\overset{_{\wedge }}{\varphi }\ \ \text{%
and \ }\overset{_{\vee }}{B}D(\psi )=\overset{_{\vee }}{\psi },
\end{equation*}%
where%
\begin{equation*}
\begin{array}{l}
\overset{_{\wedge }}{\varphi }=\xi _{u}+\xi _{u}\wedge \xi _{u}+\cdots +\xi
_{u}^{n}+\cdots \\ 
\\ 
{\xi }=\varphi +\varphi \wedge _{\ell }\varphi +\cdots +\varphi ^{n}+\cdots%
\end{array}%
\text{ \ and\ \ }%
\begin{array}{l}
\overset{_{\vee }}{\psi }=\zeta _{v}+\zeta _{v}\vee \zeta _{v}+\cdots +\zeta
_{v}^{n}+\cdots \\ 
\\ 
{\zeta }=\psi +\psi \vee _{\ell }\psi +\cdots +\psi ^{n}+\cdots .%
\end{array}%
\end{equation*}%
Then define%
\begin{equation*}
BD:C^{\ast }(P;\mathbf{U}_{+})\times C^{\ast }(P;\mathbf{U}_{+})\rightarrow
C^{\ast }(P;\mathbf{U}_{+})\times C^{\ast }(P;\mathbf{U}_{+})
\end{equation*}%
on a pair $\varphi \times \psi $ by%
\begin{equation*}
BD(\varphi \times \psi )=(\overset{_{\wedge }}{\pi }\circ \overset{_{\wedge }%
}{B}D)(\varphi ^{c}+\tau \psi )\times (\overset{_{\vee }}{\pi }\circ \overset%
{_{\vee }}{B}D)(\psi ^{a}+\tau \varphi ).\vspace{0.1in}
\end{equation*}

\begin{theorem}
\label{biderivative}Given $\dsum\limits_{\left( m,n\right) \in \mathbb{N}%
^{2}\smallsetminus \mathbf{1}}\theta _{n,m}\in {U,}$ \ there is a unique
fixed point 
\begin{equation}
\varphi \times \psi ={B}D(\varphi \times \psi )  \label{fixedpoint}
\end{equation}%
such that 
\begin{equation}
\begin{array}{rlll}
\varphi _{u_{0}}(e_{m}) & = & \theta _{1,m}, & m\geq 2 \\ 
\varphi _{u\cap v}(e_{m}) & = & \sum_{n\geq 2}\theta _{n,m}, & m\geq 2 \\ 
\psi _{v_{0}}(e_{n}) & = & \theta _{n,1}, & n\geq 2 \\ 
\psi _{u\cap v}(e_{n}) & = & \sum_{m\geq 2}\theta _{n,m}, & n\geq 2.%
\end{array}
\label{initial}
\end{equation}
\end{theorem}

Before proving this theorem, we remark that the existence of a fixed point $%
\varphi \times \psi $ for $BD$ is a deep generalization of the following
classical fact: If a map $h$ is (co)multiplicative (or a (co)derivation),
restricting $h$ to generators and (co)extend- ing as a (co)algebra map (or
as a (co)derivation) recovers $h.$ These classical (co)multi- plicative or $%
(f,g)$-(co)derivation extension procedures appear here as restrictions (\ref%
{fixedpoint}) to $P_{1}$ (a point) or to $P_{2}$ (an interval). Restricting (%
\ref{fixedpoint}) to a general permutahedron $P_{n}$ gives a new extension
procedure whose connection with the classical ones is maintained by the
compatibility of the canonical cellular projection $P_{n}\rightarrow I^{n-1}$
with diagonals. Let us proceed with a proof of Theorem \ref{biderivative}.%
\vspace{0.1in}

\begin{proof}
Define $BD^{(1)}=BD$ and $BD^{(n+1)}=BD\circ BD^{(n)},\ n\geq 1.$ Let 
\begin{equation*}
\overset{_{\wedge }}{F}_{n}\mathbf{U}=\bigoplus_{n<|\mathbf{x}|}\mathbf{U}_{%
\mathbf{x}}^{\mathbf{y}}\ \text{and}\ \overset{_{\vee }}{F}_{n}\mathbf{U}%
=\bigoplus_{n<|\mathbf{y}|}\mathbf{U}_{\mathbf{x}}^{\mathbf{y}}.
\end{equation*}%
A straightforward check shows that for each $n\geq 1,$ 
\begin{equation*}
BD^{(n+1)}=BD^{(n)}\text{ modulo }\overset{_{\wedge }}{F}_{n}C^{\ast }(P;{%
\mathbf{U}})\times \overset{_{\vee }}{F}_{n}C^{\ast }(P;{\mathbf{U}}).
\end{equation*}%
So define 
\begin{equation*}
{D}=\lim\limits_{\longrightarrow }BD^{(n)}.
\end{equation*}%
Clearly, $BD\circ {}D=D.$

Let $\varphi _{u\cap v}\in \overset{_{\wedge }}{T}op(\mathbf{U}_{+})$ and $%
\psi _{u\cap v}\in \overset{_{\vee }}{T}op(\mathbf{U}_{+})$ be the two
cochains uniquely defined by (\ref{initial}) and supported on the
appropriate faces. Then 
\begin{equation*}
\varphi \times \psi =D(\left( \varphi ^{c}+\varphi _{u\cap v}\right) \times
\left( \psi ^{a}+\psi _{u\cap v}\right) )
\end{equation*}%
is the (unique) solution of (\ref{fixedpoint}).\vspace{0.1in}
\end{proof}

\subsection{The biderivative operator on $\mathbf{U}$}

Let $\widetilde{Bd}:\mathbf{U}_{+}\times \mathbf{U}_{+}\rightarrow \mathbf{U}%
_{+}\times \mathbf{U}_{+}$ be the operator given by the composition 
\begin{equation*}
\begin{array}{ccc}
\mathbf{U}_{+}\times \mathbf{U}_{+} & \overset{\widetilde{Bd}}{%
\longrightarrow } & \mathbf{U}_{+}\times \mathbf{U}_{+} \\ 
\parallel &  & \parallel \\ 
\overset{_{\wedge }}{T}op(\mathbf{U}_{+})\times \overset{_{\vee }}{T}op(%
\mathbf{U}_{+}) & \underset{BD}{\longrightarrow } & \overset{_{\wedge }}{T}%
op(\mathbf{U}_{+})\times \overset{_{\vee }}{T}op(\mathbf{U}_{+}),%
\end{array}%
\end{equation*}%
where the vertical maps are canonical identification bijections and $BD$ is
its restriction to $\overset{_{\wedge }}{T}op(\mathbf{U}_{+})\times \overset{%
_{\vee }}{T}op(\mathbf{U}_{+})$. For $A\in \mathbf{U}_{+},$ let $A_{1}\times
A_{2}=\widetilde{Bd}\left( A\times A\right) $ and define operators $\overset{%
_{\wedge }}{B}d,\overset{_{\vee }}{B}d:\mathbf{U}_{+}\rightarrow \mathbf{U}%
_{+}$ by 
\begin{equation*}
\overset{_{\wedge }}{B}d(A)=A_{1}\text{ \ and \ }\overset{\vee }{B}%
d(A)=A_{2}.
\end{equation*}

Given an operator $F:\mathbf{U}\rightarrow \mathbf{U}$ and a submodule $%
\mathbf{U}_{\epsilon }\subset \mathbf{U,}$ denote the composition of $F$
with the projection $\mathbf{U}\rightarrow \mathbf{U}_{\epsilon }$ by $%
F_{\epsilon }$. Define the operator $Bd_{+}:\mathbf{U}_{+}\rightarrow 
\mathbf{U}_{+}$ as the sum 
\begin{equation*}
Bd_{+}=Id_{u\cap v}+\overset{_{\wedge }}{B}d_{u_{0}\oplus v}+\overset{_{\vee
}}{B}d_{u\oplus v_{0}}.
\end{equation*}%
Note that $\overset{_{\wedge }}{B}d_{u_{0}}\left( \theta \right) $ is the
cofree coextension of $\theta \in U_{1,\ast }$ as a coderivation of $T^{c}H;$
dually, $\overset{_{\vee }}{B}d_{v_{0}}\left( \eta \right) $ is the free
extension of $\eta \in U_{\ast ,1}$ as a derivation of $T^{a}H.$

On the other hand, observe that $U\cap \mathbf{U}_{0}=U_{1,1}.$ Given $A\in
U_{1,1},$ $1\leq i\leq q$ and $1\leq j\leq p,$ let ${A}_{ij}^{q\times
p}=\left( a_{k\ell }\right) \in \mathbf{U}_{\mathbf{1}^{p}}^{\mathbf{1}^{q}}$
be the $q\times p$ monomial such that%
\begin{equation*}
a_{k\ell }=\left\{ 
\begin{array}{ll}
A, & \text{if}\ (k,\ell )=(i,j), \\ 
Id, & \text{otherwise.}%
\end{array}%
\right.
\end{equation*}%
Define $Bd_{0}:U_{1,1}\rightarrow TTU_{1,1}$ by%
\begin{equation*}
Bd_{0}(A)=\sum\limits_{\substack{ 1\leq i\leq q,\text{ }1\leq j\leq p  \\ %
p,q\geq 1}}{A}_{ij}^{q\times p}.
\end{equation*}%
Then $Bd_{0}\left( A\right) $ is the free linear extension of $A$ as a
(co)derivation of $TTH$.

We establish the following fundamental notion:

\begin{definition}
The \underline{biderivative operator} 
\begin{equation*}
Bd:\mathbf{U}\rightarrow \mathbf{U}
\end{equation*}%
associated with the universal PROP $U$ is the sum 
\begin{equation*}
Bd=Bd_{0}+Bd_{+}:\mathbf{U}_{0}\oplus \mathbf{U}_{+}\rightarrow \mathbf{U}%
_{0}\oplus \mathbf{U}_{+}.
\end{equation*}%
An element $A\in \mathbf{U}$ is a \underline{biderivative} if $A=Bd(A).$
\end{definition}

\noindent Restating Theorem \ref{biderivative} in these terms we have:

\begin{proposition}
Every element $\omega={\sum\nolimits_{i,j\geq1}}\omega_{j,i}\in U$ has a
unique biderivative $d_{\omega}\in{TTU}.$
\end{proposition}

\noindent Thus the biderivative can be viewed as a non-linear map $%
d_{-}:U\rightarrow TTU.$

\subsection{The $\circledcirc $-product on $U$}

The biderivative operator allows us to extend Gerstenhaber's (co)operation 
\cite{Gersten} $\circ :U_{\ast ,1}\oplus {U_{1,\ast }\rightarrow U}$ to a
(non-bilinear) operation 
\begin{equation}
\circledcirc :U\times U\rightarrow U  \label{operation1}
\end{equation}%
defined for $\theta \times \eta \in U\times U$ by the composition 
\begin{equation*}
\circledcirc :U\times U\overset{d_{\theta }\times d_{\eta }}{\longrightarrow 
}\mathbf{U}\times \mathbf{U}\overset{\Upsilon }{\longrightarrow }\mathbf{U}%
\overset{pr}{\rightarrow }U,
\end{equation*}%
where the last map is the canonical projection. The following is now obvious:

\begin{proposition}
The $\circledcirc$ operation (\ref{operation1}) acts bilinearly only on the
submodule $U_{\ast,1}\oplus{U_{1,\ast}}.$
\end{proposition}

\begin{remark}
The bilinear part of the $\circledcirc $ operation, i.e., its restriction to 
$U_{\ast ,1}\oplus {U_{1,\ast },}$ is completely determined by the
associahedra $K$ (rather than the permutahedra) and induces the cellular
projection $P_{n}\rightarrow K_{n+1}$ due to A. Tonks \cite{Tonks}.
\end{remark}

\begin{example}
\label{example4}Throughout this example the symbol \textquotedblleft $\,1$%
\textquotedblright\ denotes the identity. Consider a graded module $H$
together with maps $d=\theta _{1,1},$ $\mu =\theta _{1,2},$ $\theta =\theta
_{2,2},$ $\Delta =\theta _{2,1}\in End\left( TH\right) .$ Let us compute the
biderivative of $\omega =d+\mu +\theta +\Delta $ and its $\circledcirc $%
-square. Consider the pair of cochains $\varphi \times \psi \in \overset{%
\wedge }{T}op\left( \mathbf{U}_{+}\right) \times \overset{\vee }{T}op\left( 
\mathbf{U}_{+}\right) $ supported on $e_{2}\times e_{2}$ such that $\varphi
\left( e_{2}\right) =\mu +\theta $ and $\psi \left( e_{2}\right) =\theta
+\Delta .$ Then 
\begin{equation*}
\begin{array}{l}
\varphi ^{c}\left( e_{2}+e_{21}+e_{12}+\cdots \right) =\mu +\left[ \mu \text{
}1\right] +\left[ 1\text{ }\mu \right] +\cdots \in \mathbf{U}_{u_{0}},%
\vspace*{0.1in} \\ 
\psi ^{a}\left( e_{2}+e_{21}+e_{12}+\cdots \right) =\Delta +\QTATOPD[ ] {%
\Delta }{1}+\QTATOPD[ ] {1}{\Delta }+\cdots \in \mathbf{U}_{v_{0}}\text{ \
and}\vspace*{0.1in} \\ 
\tau \left( \varphi \right) \left( e_{2}\right) =\varphi _{u\cap v}\left(
e_{2}\right) =\theta =\psi _{u\cap v}\left( e_{2}\right) =\tau \left( \psi
\right) \left( e_{2}\right) 
\end{array}%
\end{equation*}%
Set $\alpha =\varphi ^{c}+\tau \psi $ and $\beta =\psi ^{a}+\tau \varphi ;$
then 
\begin{eqnarray*}
\left( \alpha \wedge _{\ell }\alpha \right) \left( C_{\ast }P\right) 
&=&\left( \mu +\theta \right) \left( \mu \otimes 1+1\otimes \mu \right)
+\cdots \text{ \ and\vspace*{0.1in} } \\
\left( \beta \vee _{\ell }\beta \right) \left( C_{\ast }P\right)  &=&\left(
\Delta \otimes 1+1\otimes \Delta \right) \left( \theta +\Delta \right)
+\cdots \vspace*{0.1in}
\end{eqnarray*}%
Furthermore, the projections $\alpha _{u}=\alpha $ and $\beta _{v}=\beta $
so that 
\begin{equation*}
\xi _{u}=\alpha +\alpha \wedge _{\ell }\alpha +\cdots \text{ \ and \ }\zeta
_{v}=\beta +\beta \vee _{\ell }\beta +\cdots .
\end{equation*}%
Then $BD\left( \varphi \times \psi \right) =\hat{\varphi}\times \check{\psi},
$ where 
\begin{equation*}
\hat{\varphi}=\xi _{u}+\xi _{u}\wedge \xi _{u}+\cdots \ \ \text{and \ }%
\check{\psi}=\zeta _{v}+\zeta _{v}\vee \zeta _{v}+\cdots .
\end{equation*}%
Now $\ \hat{\varphi}_{u_{0}\oplus v}\left( C_{\ast }P\right) =\mu +\mu
\left( \mu \otimes 1+1\otimes \mu \right) +\theta +\QTATOPD[ ] {\theta }{%
\theta }+\QTATOPD[ ] {\theta }{\mu }+\QTATOPD[ ] {\mu }{\theta }+\QTATOPD[ ]
{\mu }{\mu }+\cdots $\ \ and$\vspace*{0.1in}$\linebreak $\check{\psi}%
_{u\oplus v_{0}}\left( C_{\ast }P\right) =\theta +\left[ \theta \text{ }%
\theta \right] +\left[ \Delta \text{ }\theta \right] +\left[ \theta \text{ }%
\Delta \right] +\left[ \Delta \text{\ }\Delta \right] +\Delta +\left( \Delta
\otimes 1+1\otimes \Delta \right) \Delta +\cdots \vspace*{0.1in}$\linebreak
so that $\ Bd_{+}\left( \omega \right) =\theta +\left( \hat{\varphi}%
_{u_{0}\oplus v}+\check{\psi}_{u\oplus v_{0}}\right) \left( C_{\ast
}P\right) .$ \ Finally, we adjoin the linear$\vspace*{0.1in}$\linebreak
extension of the differential $d$ in $TTU_{1,1}$ and obtain 
\begin{eqnarray*}
d_{\omega } &=&d+\left[ d\text{ }1\right] +\left[ 1\text{ }d\right] +\cdots +%
\QTATOPD[ ] {d}{1}+\QTATOPD[ ] {1}{d}+\cdots +\mu +\Delta +\cdots  \\
&&+\theta +\mu \left( \mu \otimes 1+1\otimes \mu \right) +\cdots +\left(
\Delta \otimes 1+1\otimes \Delta \right) \Delta +\cdots \vspace*{0.05in} \\
&&+\QTATOPD[ ] {\theta }{\theta }+\QTATOPD[ ] {\theta }{\mu }+\QTATOPD[ ] {%
\mu }{\theta }+\QTATOPD[ ] {\mu }{\mu }+\cdots +\left[ \theta \text{ }\theta %
\right] +\left[ \Delta \text{ }\theta \right] +\left[ \theta \text{ }\Delta %
\right] +\left[ \Delta \text{\ }\Delta \right] +\cdots .
\end{eqnarray*}%
Then (up to sign),\vspace*{-0.02in} 
\begin{eqnarray*}
\omega \circledcirc \omega  &=&\left( \QTATOPD[ ] {d}{1}+\QTATOPD[ ] {1}{d}%
\right) \cdot \theta +\theta \cdot \left( \left[ d\text{ }1\right] +\left[ 1%
\text{\ }d\right] \right) +\Delta \cdot \mu +\QTATOPD[ ] {\mu }{\mu }\cdot %
\left[ \Delta \text{\ }\Delta \right] +\hspace*{0.3in}\vspace*{0.05in} \\
&&+\QTATOPD[ ] {\mu }{\mu }\cdot \left( \left[ \Delta \text{ }\theta \right]
+\left[ \theta \text{ }\Delta \right] \right) +\theta \cdot \left( \left[ 1%
\text{ }\mu \right] +\left[ \mu \text{ }1\right] \right) +\vspace*{0.05in} \\
&&+\left( \QTATOPD[ ] {\theta }{\mu }+\QTATOPD[ ] {\mu }{\theta }\right)
\cdot \left[ \Delta \text{\ }\Delta \right] +\left( \QTATOPD[ ] {\Delta }{1}+%
\QTATOPD[ ] {1}{\Delta }\right) \cdot \theta +\cdots .
\end{eqnarray*}%
Some low dimensional relations implied by $\omega \circledcirc \omega =0$
are (up to sign):%
\begin{eqnarray*}
\left( d\otimes 1+1\otimes d\right) \theta +\theta \left( d\otimes
1+1\otimes d\right)  &=&\Delta \mu -\left( \mu \otimes \mu \right) \sigma
_{2,2}\left( \Delta \otimes \Delta \right) \vspace*{0.1in} \\
\left( \mu \otimes \mu \right) \sigma _{2,2}\left( \Delta \otimes \theta
+\theta \otimes \Delta \right)  &=&\theta \left( \mu \otimes 1+1\otimes \mu
\right) \vspace*{0.1in} \\
\left( \mu \otimes \theta +\theta \otimes \mu \right) \sigma _{2,2}\left(
\Delta \otimes \Delta \right)  &=&\left( \Delta \otimes 1+1\otimes \Delta
\right) \theta .\vspace*{0.1in}
\end{eqnarray*}%
In fact, if $\omega \circledcirc \omega =0$ then $\left( H,\omega \right) $
is an \textquotedblleft $A_{\infty }$-bialgebra.\textquotedblright 
\end{example}

\section{$A_{\infty }$-bialgebras}

In this section we define the notion of an $A_{\infty }$-bialgebra. Our
approach extends the definition of an $A_{\infty }$-(co)algebra in terms of
Gerstenhaber's (co)operation. Roughly speaking, an $A_{\infty }$-bialgebra
is a graded $R$-module $H$ equipped with compatible $A_{\infty }$-algebra
and $A_{\infty }$-coalgebra structures. Structural compatibility of the
operations in an $A_{\infty }$-bialgebra is determined by the $\circledcirc $
operation (\ref{operation1}). Before stating the definition, we mention
three natural settings in which $A_{\infty }$-bialgebras appear (details
appear in the sequel \cite{SU3}).\vspace*{0.1in}

\noindent(1) Let $X$ be a space and let $C_{\ast}(X)$ denote the simplicial
singular chain complex of $X.$ Although Adams' cobar construction $\Omega
C_{\ast}(X)$ is a (strictly coassociative) DG Hopf algebra \cite{Baues1}, 
\cite{CM}, \cite{KS1}, it seems impossible to introduce a strictly
coassociative coproduct on the double cobar construction $\Omega^{2}C_{\ast
}(X)$. Instead there is an $A_{\infty}$-coalgebra structure on $\Omega
^{2}C_{\ast}(X)$ that is compatible with the product and endows $\Omega
^{2}C_{\ast}(X)$ with an $A_{\infty}$-bialgebra structure.\vspace{0.1in}

\noindent(2) If $H$ is a graded bialgebra and $\rho:RH\longrightarrow H$ is
a (bigraded) multiplicative resolution, it is difficult to introduce a
strictly coassociative coproduct on $RH$ in such a way that $\rho$ is a map
of bialgebras. However, there exists an $A_{\infty}$-bialgebra structure on $%
RH$ such that $\rho$ is a morphism of $A_{\infty}$-bialgebras.\vspace{0.1in}

\noindent(3) If $A$ is any DG bialgebra, its homology $H(A)$ has a canonical 
$A_{\infty}$-bialgebra structure.\vspace*{0.1in}

The definition of an $A_{\infty }$-bialgebra $H$ uses the $\circledcirc $%
-operation on $U_{H}$ to mimic the definition of an $A_{\infty }$-algebra.

\begin{definition}
An \underline{$A_{\infty }$-bialgebra} is a graded $R$-module $H$ equipped
with operations 
\begin{equation*}
\{\omega ^{j,i}\in Hom^{i+j-3}(H^{\otimes i},H^{\otimes j})\}_{i,j\geq 1}
\end{equation*}%
such that $\omega =\tsum_{i,j\geq 1}\omega ^{j,i}\in U$ satisfies $\omega
\circledcirc \omega =0.$
\end{definition}

Here are some of the first structural relations among the operations in an $%
A_{\infty}$-bialgebra:

\begin{equation*}
\begin{array}{llll}
d\omega ^{2,2}= & \omega ^{2,1}\omega ^{1,2}-\left( \omega ^{1,2}\otimes
\omega ^{1,2}\right) \sigma _{2,2}(\omega ^{2,1}\otimes \omega ^{2,1}) &  & 
\\ 
&  &  &  \\ 
d\omega ^{3,2}= & \omega ^{3,1}\omega ^{1,2}+(\omega ^{2,1}\otimes
1-1\otimes \omega ^{2,1})\omega ^{2,2} &  &  \\ 
& -(\omega ^{1,2}\otimes \omega ^{1,2}\otimes \omega ^{1,2})\sigma _{3,2} 
\left[ \omega ^{3,1}\otimes (1\otimes \omega ^{2,1})\omega ^{2,1}+(\omega
^{2,1}\otimes 1)\omega ^{2,1}\otimes \omega ^{3,1}\right] &  &  \\ 
& +\left[ (\omega ^{2,2}\otimes \omega ^{1,2}-\omega ^{1,2}\otimes \omega
^{2,2})\right] \sigma _{2,2}(\omega ^{2,1}\otimes \omega ^{2,1}) &  &  \\ 
&  &  &  \\ 
d\omega ^{2,3}= & -\omega ^{2,1}\omega ^{1,3}+\omega ^{2,2}(1\otimes \omega
^{1,2}-\omega ^{1,2}\otimes 1) &  &  \\ 
& +[\omega ^{1,3}\otimes \omega ^{1,3}(1\otimes \omega ^{1,3})+\omega
^{1,3}(\omega ^{1,2}\otimes 1)\otimes \omega ^{1,3}]\sigma _{2,3}(\omega
^{2,1}\otimes \omega ^{2,1}\otimes \omega ^{2,1}) &  &  \\ 
& +(\omega ^{1,2}\otimes \omega ^{1,2})\sigma _{2,2}(\omega ^{2,1}\otimes
\omega ^{2,2}-\omega ^{2,2}\otimes \omega ^{2,1}). &  & 
\end{array}%
\end{equation*}

\begin{example}
The structure of an $A_{\infty }$-bialgebra whose initial data consists of a
strictly coassociative coproduct $\Delta :H\rightarrow H^{\otimes 2}$
together with $A_{\infty }$-algebra operations $m_{i}:H^{\otimes
i}\rightarrow H,$ $i\geq 2,$ is determined as in Example \ref{example4} but
with $\varphi \left( e_{i}\right) =m_{i},$ $\psi \left( e_{2}\right) =\Delta
.$ This time the action of $\tau $ is trivial since all initial maps lie in $%
\mathbf{U}_{u_{0}\oplus v_{0}}$ and we obtain the following structure
relation for each $i\geq 2$: 
\begin{equation*}
\left( \xi _{u}\wedge \xi _{u}\right) \left( e_{i}\right) \cdot \lbrack 
\underset{i\text{ factors}}{\underset{}{\Delta \text{ }\cdots \text{ }\Delta 
}}]=\Delta \cdot m_{i}.
\end{equation*}%
Indeed, the classical bialgebra relation appears when $i=2.$
\end{example}

We conclude with a statement of our main theorem (the definition of an $%
A_{\infty }$-bialgebra morphism appears in the sequel \cite{SU3}).

\begin{theorem}
Let $A$ be an $A_{\infty }$-bialgebra; if the ground ring $R$ is not a
field, assume that the homology $H=H(A)$ is torsion-free. Then $H$ inherits
a canonical bialgebra structure that extends to an $A_{\infty }$-bialgebra
structure $\{\omega ^{j,i}\}_{i,j\geq 1}\ $with $\omega ^{1,1}=0.$
Furthermore, there is a map of $A_{\infty }$-bialgebras%
\begin{equation*}
\digamma =\{\digamma ^{j,i}\}_{i,j\geq 1}:H\Longrightarrow A,
\end{equation*}%
with $\digamma ^{j,i}\in Hom^{i+j-2}(H^{\otimes i},A^{\otimes j}),$ such
that $\digamma ^{1,1}:H\rightarrow A$ is a map of DGM's inducing an
isomorphism on homology.
\end{theorem}


\begin{thebibliography}{99}
\bibitem{Baues1} H. J. Baues, The cobar construction as a Hopf algebra and
the Lie differential, \textit{Invent. Math.} \textbf{132} (1998) 467-489.

\bibitem{CM} G. Carlsson and R. J. Milgram, Stable homotopy and iterated
loop spaces, \textit{Handbook of Algebraic Topology (I. M. James, ed.)},
North-Holland (1995), 505-583.

\bibitem{Gersten} M. Gerstenhaber and S. D. Schack, Algebras, bialgebras,
quantum groups, and algebraic deformations, \textit{Contemporary Math.} 
\textbf{134,} A. M. S., Providence (1992), 51-92.

\bibitem{Gugenheim} V.K.A.M. Gugenheim, On a perturbation theory for the
homology of the loop space, \textit{J. Pure Appl. Algebra,} \textbf{25}
(1982), 197-205.

\bibitem{Kadeishvili1} T. Kadeishvili, On the homology theory of fibre
spaces, \textit{Russian Math. Survey}, \textbf{35} (1980), 131-138.

\bibitem{Kimura1} T. Kimura, J. Stasheff and A. Voronov, On operad
structures of moduli spaces and string theory, \textit{Comm. Math. Physics}, 
\textbf{171} (1995), 1-25.

\bibitem{Kimura2} T. Kimura, A. Voronov and G. Zuckerman, Homotopy
Gerstenhaber algebras and topological field theory, \textit{Operads:
Proceedings of Renaissance Conferences (J.-L. Loday, J. Stasheff and A.
Voronov, eds.)}, A. M. S. Contemp. Math. \textbf{202} (1997), 305-334.

\bibitem{KS1} T. Kadeishvili and S. Saneblidze, A cubical model of a
fibration, \textit{J. Pure Appl. Algebra}, \textbf{196} (2005), 203-228.

\bibitem{Lada} T. Lada and M. Markl, Strongly homotopy Lie algebras, \textit{%
Communications in Algebra,} \textbf{23} (1995), 2147-2161.

\bibitem{Lin} J. P. Lin, $H$-spaces with finiteness conditions, \textit{%
Handbook of Algebraic Topology\ (I. M. James, ed.)}, North Holland,
Amsterdam (1995), 1095-1141.

\bibitem{Loday} J.-L. Loday and M. Ronco, Hopf algebra of the planar binary
trees, \textit{Adv. in Math. }\textbf{139, }No. 2 (1998), 293-309.

\bibitem{Markl} M. Markl, A cohomology theory for $A(m)$-algebras and
applications, \textit{J. Pure and Appl. Algebra,} \textbf{83} (1992),
141-175.

\bibitem{Penkava1} M. Penkava and A. Schwarz, On some algebraic structures
arising in string theory, \textquotedblleft Conf. Proc. Lecture Notes Math.
Phys., III,\ Perspectives in Math. Physics\textit{,\textquotedblright }
International Press, Cambridge (1994), 219-227.

\bibitem{Penkava2} --------------------, $A_{\infty }$ algebras and the
cohomology of moduli spaces, \textquotedblleft Lie Groups and Lie Algebras:
E. B. Dykin's Seminar,\textquotedblright\ A. M. S. Transl. Ser. 2 \textbf{169%
} (1995), 91-107.

\bibitem{Saneblidze} S. Saneblidze, On the homotopy classification of spaces
by the fixed loop space homology, \textit{Proc. A. Razmadze Math. Inst.}, 
\textbf{119} (1999), 155-164.

\bibitem{SU1} S. Saneblidze and R. Umble, A Diagonal on the associahedra,
preprint AT/0011065, November 2000.

\bibitem{SU2} --------------, Diagonals on the permutahedra, multiplihedra
and associahedra, \textit{J. Homology, Homotopy and Appl}., \textbf{6 }(1)
(2004), 363-411.

\bibitem{SU3} --------------, Matrons and the category of $A_{\infty }$%
-bialgebras, in preparation.

\bibitem{Shnider} S. Shnider and S. Sternberg, \textquotedblleft Quantum
Groups: From Coalgebras to Drinfeld Algebras,\textquotedblright\
International Press, Boston (1993).

\bibitem{borya} B. Shoikhet, The CROCs, non-commutative deformations, and
(co)associative bialgebras, preprint QA/0306143.

\bibitem{Smith} J. R. Smith, \textquotedblleft Iterating the Cobar
Construction,\textquotedblright\ Memoirs of the A. M. S. \textbf{109},
Number 524 (1994).

\bibitem{Stasheff} J. D. Stasheff, Homotopy associativity of $H$-spaces I,
II, \textit{Trans. A. M. S.} \textbf{108} (1963), 275-312.

\bibitem{Stasheff2} ------------------, \textquotedblleft $H$-spaces from a
Homotopy Point of View,\textquotedblright\ SLNM 161, Springer, Berlin (1970).

\bibitem{Tonks} A. Tonks, Relating the associahedron and the permutohedron,
\textquotedblleft Operads: Proceedings of the Renaissance Conferences
(Hartford CT / Luminy Fr 1995),\textquotedblright\ Contemporary Mathematics 
\textbf{202} (1997), pp.33-36 .

\bibitem{Umble} R. N. Umble, The deformation complex for differential graded
Hopf algebras, \textit{J. Pure Appl. Algebra,} \textbf{106} (1996), 199-222.

\bibitem{Umble2} --------------------, In Search of higher homotopy Hopf
algebras, lecture notes.

\bibitem{Wiesbrock} H. W. Wiesbrock, A Note on the construction of the $%
C^{\ast }$-Algebra of bosonic strings, \textit{J. Math. Phys.} \textbf{33}
(1992), 1837-1840.

\bibitem{Zwiebach} B. Zwiebach, Closed string field theory; quantum action
and the Batalin-Vilkovisky master equation, \textit{Nucl. Phys. B}. \textbf{%
390} (1993), 33-152.
\end{thebibliography}
\end{document}